\documentclass[11pt]{article}
\usepackage{graphicx,subfigure,ifthen,color,algorithm,palatino}
\usepackage{amsmath,amsthm,amssymb,amsfonts,verbatim}
\usepackage{mathrsfs,accents,setspace,booktabs,subfig}
\usepackage{soul,cancel}
\DeclareMathAlphabet\mathbfcal{OMS}{cmsy}{b}{n}
\newboolean{UnBlinded}
\newboolean{DoubleSpaced}
\setboolean{UnBlinded}{true}
\setboolean{DoubleSpaced}{false}
\def\myGeneralized{Generalised}
\def\mygeneralized{generalised}
\def\myOP{O_p}
\def\myoP{o_p}
\def\thick#1{\hbox{\rlap{$#1$}\kern0.25pt\rlap{$#1$}\kern0.25pt$#1$}}
\def\bone{\boldsymbol{1}}
\def\bI{\boldsymbol{I}}
\def\bQ{\boldsymbol{Q}}
\def\bS{\boldsymbol{S}}
\def\bX{\boldsymbol{X}}
\def\bY{\boldsymbol{Y}}

\def\Csc{{\mathcal C}}
\def\convprob{\stackrel{P}{\to}}
\def\Var{\mbox{Var}}
\def\argmindum{\mathop{\mbox{\rm argmin}}}
\def\argmin#1{\argmindum_{#1}}
\def\bib{\vskip12pt\par\noindent\hangindent=1 true cm\hangafter=1}
\def\ba{\boldsymbol{a}}
\def\bb{\boldsymbol{b}}

\def\be{\boldsymbol{e}}

\def\bg{\boldsymbol{g}}
\def\bh{\boldsymbol{h}}
\def\phihat{{\widehat\phi}}
\def\bSigmahat{{\widehat\bSigma}}
\def\bPhihat{{\widehat\bPhi}}
\def\bv{\boldsymbol{v}}

\def\bx{\boldsymbol{x}}
\def\by{\boldsymbol{y}}

\def\tcbk#1{\textcolor{black}{#1}}
\def\tr{\mbox{tr}}
\def\Psc{{\mathcal P}}
\def\bU{\boldsymbol{U}}
\def\bzero{\boldsymbol{0}}
\def\bSigma{\boldsymbol{\Sigma}}
\def\bu{\boldsymbol{u}}
\def\argmaxdum{\mathop{\mbox{\rm argmax}}}
\def\argmax#1{\argmaxdum_{#1}}
\def\vech{\mbox{\rm vech}}
\def\Cov{\mbox{\rm Cov}}
\def\bLambda{\boldsymbol{\Lambda}}
\def\sumim{\sum_{i=1}^m}
\def\vecof{\mbox{\rm vec}}
\def\Gsc{{\mathcal G}}
\def\Hsc{{\mathcal H}}
\def\bbetahat{{\widehat\bbeta}}
\def\smbbeta{{\thick{\scriptstyle{\beta}}}}
\def\smbSigma{{\thick{\scriptstyle{\Sigma}}}}
\def\bPsi{\boldsymbol{\Psi}}
\def\bPhi{\boldsymbol{\Phi}}
\def\bOmega{\boldsymbol{\Omega}}
\def\bA{\boldsymbol{A}}
\def\bB{\boldsymbol{B}}
\def\bC{\boldsymbol{C}}
\def\bD{\boldsymbol{D}}
\def\bK{\boldsymbol{K}}
\def\bM{\boldsymbol{M}}
\def\bO{\boldsymbol{O}}
\def\betahat{{\widehat\beta}}
\def\smhalf{{\textstyle{\frac{1}{2}}}}

\def\VECTORone{\psi_1(\bU)}
\def\VECTORtwo{\psi_2(\bU)}
\def\VECTORthree{\psi_3(\bU)}
\def\VECTORfour{\psi_4(\bU)}
\def\MATRIXfive{\bPsi_5(\bU)}
\def\MATRIXsix{\bPsi_6(\bU)}
\def\MATRIXseven{\bPsi_7(\bU)}
\def\MATRIXeight{\bPsi_8(\bU)}
\def\MATRIXnine{\bPsi_9(\bU)}
\def\bOmegaAA{\bOmega_{\mbox{\tiny AA}}}
\def\bOmegaAB{\bOmega_{\mbox{\tiny AB}}}
\def\bOmegaBB{\bOmega_{\mbox{\tiny BB}}}
\def\bOmegadAAA{\bOmega'_{\mbox{\tiny AAA}}}
\def\bOmegadAAB{\bOmega'_{\mbox{\tiny AAB}}}
\def\bLambdaAA{\bLambda_{\mbox{\tiny AA}}}
\def\bLambdaAB{\bLambda_{\mbox{\tiny AB}}}
\def\bOmegaAAhat{{\widehat\bOmega}_{\mbox{\tiny AA}}}
\def\bOmegaABhat{{\widehat\bOmega}_{\mbox{\tiny AB}}}
\def\bOmegaBBhat{{\widehat\bOmega}_{\mbox{\tiny BB}}}
\def\AsyCovHat{\widehat{\mbox{Asy.Cov}}}
\def\bLambdaAAhat{{\widehat{\bLambda}}_{\mbox{\tiny AA}}}
\def\bLambdaABhat{{\widehat{\bLambda}}_{\mbox{\tiny AB}}}
\def\bAsc{\mathbfcal{A}}
\def\bBsc{\mathbfcal{B}}
\def\bKsc{\mathbfcal{K}}
\def\bXsc{\mathbfcal{X}}
\def\bKscAA{\bKsc_{\mbox{\tiny AA}}}
\def\bKscAB{\bKsc_{\mbox{\tiny AB}}}
\def\bKscBB{\bKsc_{\mbox{\tiny BB}}}
\def\bKscBC{\bKsc_{\mbox{\tiny BC}}}

\def\bKscCC{\bKsc_{\mbox{\tiny CC}}}
\def\HscrdAAAi{\Hsc'_{\mbox{\scriptsize AAA}i}}
\def\HscrdAABi{\Hsc'_{\mbox{\scriptsize AAB}i}}
\def\dB{d_{\mbox{\tiny B}}}
\def\bX{\boldsymbol{X}}
\def\convprobBka{\stackrel{p}{\to}}
\def\bbeta{\boldsymbol{\beta}}
\def\intdR{\int_{\real^{\dR}}}

\def\bbetaMLE{{\widehat\bbeta}}
\def\betaAMLE{{\widehat\beta}_{\mbox{\scriptsize A}}}
\def\bbetaAMLE{{\widehat\bbeta}_{\mbox{\scriptsize A}}}
\def\bbetaBMLE{{\widehat\bbeta}_{\mbox{\scriptsize B}}}

\def\bSigmaMLE{{\widehat\bSigma}}

\def\buzero{\left[\begin{array}{c}\bu\\\bzero\end{array}\right]}
\def\dF{d_{\mbox{\tiny F}}}
\def\dR{d_{\mbox{\tiny R}}}
\def\intdR{\int_{\real^{\dR}}}
\def\bXA{\bX_{\mbox{\scriptsize A}}}
\def\bXB{\bX_{\mbox{\scriptsize B}}}
\def\bXij{\bX_{ij}}
\def\bXAij{\bX_{\mbox{\scriptsize A}ij}}
\def\bXBij{\bX_{\mbox{\scriptsize B}ij}}
\def\bbetaA{\bbeta_{\mbox{\scriptsize A}}}
\def\bbetaB{\bbeta_{\mbox{\scriptsize B}}}

\def\tinybbetaA{\bbeta_{\mbox{\tiny A}}}
\def\tinybbetaB{\bbeta_{\mbox{\tiny B}}}

\def\tinyvechbSigma{\vech(\bSigma)}
\def\bbetaAzero{\bbetaA^0}
\def\bbetaBzero{\bbetaB^0}
\def\bSigmaZero{\bSigma^0}
\def\bbetaZero{\bbeta^0}

\def\bSAi{\bS_{\mbox{\scriptsize A}i}}
\def\bSBi{\bS_{\mbox{\scriptsize B}i}}
\def\bSCi{\bS_{\mbox{\scriptsize C}i}}
\def\GscrAi{\Gsc_{\mbox{\scriptsize A}i}}
\def\GscrBi{\Gsc_{\mbox{\scriptsize B}i}}
\def\HscrAAi{\Hsc_{\mbox{\scriptsize AA}i}}
\def\HscrABi{\Hsc_{\mbox{\scriptsize AB}i}}
\def\HscrBBi{\Hsc_{\mbox{\scriptsize BB}i}}

\def\nablabbetaA{\nabla_{\mbox{\tiny $\tinybbetaA$}}}
\def\nablabbetaB{\nabla_{\mbox{\tiny $\tinybbetaB$}}}

\def\nablavechbSigma{\nabla_{\mbox{\scriptsize $\tinyvechbSigma$}}}

\def\linXABUi{(\bbetaA+\bU_i)^T\bXAij+\bbetaB^T\bXBij}
\def\linXABu{(\bbetaA+\bu)^T\bXAij+\bbetaB^T\bXBij}
\def\linXABUistar{(\bbetaA+\bU_i^*)^T\bXAij+\bbetaB^T\bXBij}

\def\linXABUzero{(\bbetaAzero+\bU)^T\bXA+(\bbetaBzero)^T\bXB}

\def\smhalfdisp{{\displaystyle\frac{1}{2}}}

\def\sumjni{\sum_{j=1}^{n_i}}

\def\dispsumjni{{\displaystyle\sumjni}}
\def\bgiA{\bg_{\scriptscriptstyle iA}}
\def\bgiB{\bg_{\scriptscriptstyle iB}}
\def\bgiC{\bg_{\scriptscriptstyle iC}}
\def\cS{c_{\scriptscriptstyle S}}

\def\smoon{{\textstyle{\frac{1}{n}}}}

\def\real{{\mathbb R}}

\def\betazZero{\beta_0^0}

\def\sigsqZero{(\sigma^2)^0}

\def\oon{\frac{1}{n}}

\def\myand{\&\ }

\def\betazZero{\beta_0^0}
\def\betaoZero{\beta_1^0}

\def\sigsqZero{(\sigma^2)^0}

\def\pYiGXi{p_{\mbox{\tiny $\bY_i|\bX_i$}}}

\def\dispsumjn{{\displaystyle\sum_{j=1}^{n_i}}}

\def\naturalNumbers{{\mathbb N}}
\def\EpsiSixHat{\widehat{E}\big\{\MATRIXsix\big\}}
\def\EpsiEightHat{\widehat{E}\big\{\MATRIXeight\big\}}
\def\EpsiNineHat{\widehat{E}\big\{\MATRIXnine\big\}}
\def\betaNoughtZero{\beta_0^0}
\def\betaOneZero{\beta_1^0}
\def\betaOneZero{\beta_1^0}
\def\betaTwoZero{\beta_2^0}
\def\betaThreeZero{\beta_3^0}
\def\betaFourZero{\beta_4^0}

\newtheorem{lemma}{\textbf{Lemma}}
\begin{document}

\ifthenelse{\boolean{DoubleSpaced}}{\setstretch{1.5}}{}
\vskip5mm
\centerline{\Large\bf Second Term Improvement to \myGeneralized}
\vskip1mm
\centerline{\Large\bf Linear Mixed Model Asymptotics}
\vskip5mm
\centerline{\normalsize\sc By Luca Maestrini$\null^1$, 
Aishwarya Bhaskaran$\null^2$ and Matt P. Wand$\null^3$}
\vskip5mm
\centerline{\textit{$\null^1$\tcbk{The} Australian National University, 
$\null^2$Macquarie University 
and 
$\null^3$University of Technology Sydney}}
\vskip6mm
\centerline{31st March, 2023}
\vskip6mm
\centerline{\large\bf Abstract}
\vskip2mm

A recent article on \mygeneralized\ linear mixed model asymptotics,
Jiang \textit{et al.} (2022), derived the rates of convergence
for the asymptotic variances of maximum likelihood estimators.
If $m$ denotes the number of groups and $n$ is the average within-group
sample size then the asymptotic variances have orders $m^{-1}$ and 
$(mn)^{-1}$, depending on the parameter. We extend this theory
to provide explicit forms of the $(mn)^{-1}$ second terms of the 
asymptotically harder-to-estimate parameters.
Improved accuracy of studentised confidence intervals is one consequence of our theory.

\vskip3mm
\noindent
\textit{Keywords:} Longitudinal data analysis, Maximum likelihood estimation, Studentisation.

\section{Introduction}\label{sec:intro}

\myGeneralized\ linear mixed models are a vehicle for regression 
analysis of grouped data with non-Gaussian responses such as counts 
and categorical labels. Until recently, the precise asymptotic behaviours of 
the conditional maximum likelihood estimators were not known \tcbk{for these models}.
Jiang \textit{et al.} (2022) derived leading term asymptotic variances and
showed them have orders $m^{-1}$ and $(mn)^{-1}$, depending on the parameter,
where $m$ is the number of groups and $n$ is the average within-group 
sample size. The main contribution of this article is to extend the asymptotic variance
and covariance approximations to terms in $(mn)^{-1}$ for \emph{all}
parameters. This constitutes \emph{second term improvement} to generalized\ linear 
mixed model asymptotics. The potential statistical payoffs are improved 
accuracy of confidential intervals, hypothesis tests, sample size calculations 
and optimal design.

The essence of generalized\ linear mixed models is the extension of general
linear models via the addition of random effects that allow for the handling
of correlations arising from repeated measures. There are numerous types of 
random effect structures. The most common is the two-level nested structure, 
corresponding to repeated measures within each of $m$ distinct groups. This 
version of \mygeneralized\ linear mixed models, with frequentist 
inference via maximum likelihood and its quasi-likelihood extension, is our 
focus here. Overviews of \mygeneralized\ linear mixed models are provided by 
books such as Jiang \myand Nguyen (2021), McCulloch \textit{et al.} (2008) 
and Stroup (2013). 

Suppose that a fixed effects parameter in a two-level \mygeneralized\ linear mixed model
is accompanied by a random effect. Jiang \textit{et al.} (2022) showed that the 
variance of its maximum likelihood estimator, conditional on the predictor data, 
is asymptotic to $C_1m^{-1}$ for some deterministic constant $C_1$ that depends
on the true model parameter values. The crux of this article is to extend the 
asymptotic variance approximation to $C_1m^{-1}+C_2(mn)^{-1}$ for an
additional deterministic constant $C_2$. We derive the explicit form of $C_2$
for two-level nested \mygeneralized\ linear mixed models for both maximum likelihood
and maximum quasi-likelihood situations. Even though, in general, $C_2$ does not
have a succinct form it is still usable in that operations such as studentisation
are straightforward and result in improvements in statistical utility.

For two-level nested mixed models, $(mn)^{-1}$ is the best possible rate of convergence 
for the asymptotic variance of the estimator of a model parameter. Such a rate is 
achieved by maximum likelihood estimators of fixed effects parameters unaccompanied
by random effects and dispersion parameters (e.g. Bhaskaran \myand Wand, 2023).
The current article closes the problem of obtaining the precise asymptotic
forms of the variances, up to terms in $(mn)^{-1}$, for estimation of \emph{all} model 
parameters.

Section \ref{sec:modelDescrip} describes the model under consideration and
corresponding maximum estimators. Our second term improvement
results are presented in Section \ref{sec:mainres}. Section \ref{sec:utility}
describes statistical utility due to the new asymptotic results. We 
present some corroborating numerical results in Section \ref{sec:numericRes}.
A supplement to this article contains derivational details.

\section{Model Description and Maximum Likelihood Estimation}\label{sec:modelDescrip}

Consider the class of two-parameter exponential family of density, or probability mass,
functions with generic form
\begin{equation}
p(y;\eta,\phi)= \exp[\left\{y\eta - b(\eta)+c(y)\right\}/\phi + d(y,\phi)]h(y)
\label{eq:expDens}
\end{equation}
where $\eta$ is the \emph{natural parameter} and $\phi>0$ is the \emph{dispersion parameter}.
Examples include the Gaussian density for which $b(x)= \smhalf x^2$, 
$c(x)=-\smhalf x^2$, $d(x_1,x_2)=-\smhalf\log(2\pi x_2)$ and $h(x)=I(x\in\real)$ and the 
Gamma density function for which $b(x)=-\log(-x)$, $c(x)=\log(x)$, 
$d(x_1,x_2)=-\log(x_1)-\log(x_2)/x_2-\log\Gamma(1/x_2)$ and 
$h(x)=I(x>0)$. Here $I(\Psc)=1$ if the condition $\Psc$ is 
true and $I(\Psc)=0$ if $\Psc$ is false. The Binomial and Poisson probability
mass functions are also special cases of (\ref{eq:expDens}) but with 
$\phi$ fixed at 1. When (\ref{eq:expDens}) is used in regression contexts
a common modelling extension for count and proportion responses, usually to 
account for overdispersion, is to remove the $\phi=1$ restriction and replace it
with $\phi>0$. In these circumstances $\{y\eta - b(\eta)+c(y)\}/\phi + d(y,\phi)$ 
is labelled a \emph{quasi-likelihood function} since it is not the 
logarithm of a probability mass function for $\phi\ne1$.
We use the more general quasi-likelihood terminology for the remainder 
of this article.

Consider, for observations of the random pairs $(\bXij,Y_{ij})$, $1\le i\le m$, $1\le j\le n_i$,
\mygeneralized\ linear mixed models of the form, 
\begin{equation}
\begin{array}{l}
Y_{ij}|\bXij,\bU_i\ \mbox{independent having quasi-likelihood function (2) with natural \tcbk{parameter}}\\[1ex]
\left(\bbetaZero+\left[\begin{array}{c}\bU_i\\ 
\bzero
\end{array}
\right]
\right)^T\bXij
\ \mbox{such that the}
\ \bU_i\ 
\mbox{are independent $N(\bzero,\bSigmaZero)$ \tcbk{random vectors.}}
\end{array}
\label{eq:theModel}
\end{equation}
The $\bXij$ are $\dF\times1$ random vectors corresponding to predictors.
The $\bU_i$ are $\dR\times1$ unobserved random effects vectors, where $\dR\le\dF$.
Under this set-up the first $\dR$ entries of the $\bXij$ are partnered by
a random effect. The remaining entries correspond to predictors that have
a fixed effect only. We assume that the $\bX_{ij}$ and $\bU_i$, for $1\le i\le m$ 
and $1\le j\le n_i$, are totally independent, with the $\bX_{ij}$ each having the
same distribution as the $\dF\times 1$ random vector $\bX$ 
and the $\bU_i$ each having the same distribution as the $\dR\times 1$ 
random vector $\bU$.

For any $\bbeta$ $(\dF\times1)$ and $\bSigma$ $(\dR\times\dR)$ that is 
symmetric and positive definite and conditional on the $\bX_{ij}$ data, 
the quasi-likelihood is
{\setlength\arraycolsep{1pt}
\begin{eqnarray*}
&&\ell(\bbeta,\bSigma)=\sum_{i=1}^m\sum_{j=1}^{n_i}[\{Y_{ij}(\bbeta^T\bXij+c(Y_{ij})\}/\phi
+d(Y_{ij},\phi)]-\frac{m}{2}\log|2\pi\bSigma|\\
&&\qquad\quad\quad\quad
+\sum_{i=1}^m\log\int_{\real^{\dR}}\exp\Bigg[\frac{1}{\phi}\sum_{j=1}^{n_i}\left\{Y_{ij}
\buzero^T\bXij
-b\left(\left(\bbeta+\buzero\right)^T\bXij\right)\right\}
-\smhalf\bu^T\bSigma^{-1}\bu\Bigg]\,d\bu.
\end{eqnarray*}
}
The maximum quasi-likelihood estimator of $(\bbetaZero,\bSigmaZero)$ is 
$$(\bbetaMLE,\bSigmaMLE)=\argmax{\bbeta,\bSigma}\ell(\bbeta,\bSigma).$$

Suppose that $\dF>\dR$ and consider the partition $\bbeta=[\bbetaA^T\ \bbetaB^T]^T$ 
of the fixed effects parameter vector, where $\bbetaA$ is $\dR\times 1$ and $\bbetaB$
is $(\dF-\dR)\times1$. The $\dF=\dR$ boundary case is such that $\bbetaB$ is null.
Also, let 
$\bXsc\equiv\{\bX_{ij}:1\le i\le m,\ 1\le j\le n_i\}$.
Theorem 1 of Jiang \textit{et al.} (2022) implies that, under
some mild conditions, the covariance matrices of $\bbetaAMLE$, 
$\bbetaBMLE$ and $\vech(\bSigmaMLE)$ have leading term behaviour
given by
\begin{equation}
\Cov\big(\bbetaAMLE|\bXsc\big)=\frac{\bSigmaZero\{1+o_p(1)\}}{m},\ \ 
\Cov\big(\bbetaBMLE|\bXsc\big)=\frac{\tcbk{\phi}\bLambda_{\bbetaB}\{1+o_p(1)\}}{mn},
\ \ \mbox{where}\ \ n\equiv \frac{1}{m}\sumim n_i,
\label{eq:Kerplunk}
\end{equation}
and
\begin{equation}
\Cov\big(\vech(\bSigmaMLE)|\bXsc\big)
=\frac{2\bD_{\dR}^{+}(\bSigmaZero\otimes\bSigmaZero)\bD_{\dR}^{+T}\{1+o_p(1)\}}{m}.
\label{eq:dropBears}
\end{equation}
Here
$\bLambda_{\bbetaB}$ is a  $(\dF-\dR)\times(\dF-\dR)$ matrix that depends
on $\bbeta$ and the $(\bX,\bU)$ distribution, $\bD_{\dR}$ is the matrix of 
zeroes and ones such that 
$\bD_{\dR}\vech(\bA)=\vecof(\bA)$ for all $\dR\times\dR$ symmetric matrices $\bA$
and $\bD_{\dR}^+=(\bD_{\dR}^T\bD_{\dR})^{-1}\bD_{\dR}^T$ is the 
Moore-Penrose inverse of $\bD_{\dR}$. The theory of Jiang \textit{et al.} (2022) 
also indicates a degree of asymptotic orthogonality between $\bbetaA$ and $\bbetaB$
in that $E\big\{(\bbetaAMLE-\bbetaAzero)(\bbetaBMLE-\bbetaBzero)^T|\bXsc\big\}$ has 
$\myOP\{(mn)^{-1}\}$ entries, which implies that the correlations between the 
entries of $\bbetaAMLE$ and $\bbetaBMLE$ are asymptotically negligible.

The leading term approximations of the variability in $\bbetaAMLE$ and $\vech(\bSigmaMLE)$,
given by (\ref{eq:Kerplunk}) and (\ref{eq:dropBears}), are somewhat
crude. Unlike the asymptotic covariance of $\bbetaBMLE$, they do not show the 
effect of the average within-group sample size $n$. In the next section we 
investigate their second term improvements.

\section{Two-Term Asymptotic Covariance Results}\label{sec:mainres}

We define the two-term asymptotic covariance matrix problem
to be the determination of the \emph{unique deterministic} matrices $\bM_{\smbbeta}$
and $\bM_{\smbSigma}$ such that
{\setlength\arraycolsep{1pt}
\begin{eqnarray*}
\Cov\big(\bbetaMLE|\bXsc\big)&=&
\frac{1}{m}\left[
\begin{array}{cc}
\bSigmaZero &\quad \bO \\[1ex]
\bO     &\quad  \bO
\end{array}
\right]
+\frac{\bM_{\smbbeta}\{1+\myoP(1)\}}{mn}
\quad\mbox{and}
\\[1ex]
\Cov\big(\vech(\bSigmaMLE)|\bXsc\big)&=&
\displaystyle{\frac{2\bD_{\dR}^+(\bSigmaZero\otimes\bSigmaZero)
\bD_{\dR}^{+T}}{m}}
+\frac{\bM_{\smbSigma}\{1+\myoP(1)\}}{mn}
\end{eqnarray*}
}
under reasonably mild conditions. 

An example for which a solution to the two-term asymptotic covariance
problem can be expressed relatively simply is the $\dF=2$, $\dR=1$ Poisson
quasi-likelihood special case of (\ref{eq:theModel}), with parameters
$$\bbeta=(\beta_0,\beta_1)\quad\mbox{and}\quad\bSigma=\sigma^2
\quad\mbox{and predictor variable}\quad
\bX=\left[
\begin{array}{c}
1\\[1ex]
X
\end{array}
\right]
$$
for a scalar random variable $X$. Define 
$$a_1\big(\beta_0,\beta_1,\sigma^2\big)\equiv
e^{\beta_0+\sigma^2/2}\big[E(X^2e^{\beta_1X})E(e^{\beta_1X})-\{E(Xe^{\beta_1X})\}^2\big]$$
and
$$a_2(\beta_1,\sigma^2)\equiv 
\frac{e^{\sigma^2}E\big(X^2e^{\beta_1X}\big)E\big(e^{\beta_1X}\big)
+\big(1-e^{\sigma^2}\big)E\{\big(Xe^{\beta_1X}\big)\}^2}{E\big(e^{\beta_1X}\big)}.
$$
Then the two-term covariance matrix of $(\betahat_0,\betahat_1)$ is
$$
\Cov\left(
\left[
\begin{array}{c}
\betahat_0\\[1ex]
\betahat_1
\end{array}
\right]
\Bigg|\bXsc
\right)=
\frac{1}{m}
\left[
\begin{array}{cc}
\sigsqZero &\ 0 \\[1ex]
0          &\ 0
\end{array}
\right]
+
\frac{\phi\{1+\myoP(1)\}}{a_1\big(\betazZero,\betaoZero,\sigsqZero\big)\,mn}
\left[
\begin{array}{cc}
a_2\big(\betaoZero,\sigsqZero\big)   &\ -E\big(Xe^{\betaoZero X}\big) \\[1ex]
-E\big(Xe^{\betaoZero X}\big)        &\ E\big(e^{\betaoZero X}\big)
\end{array}
\right].
$$
In other words, for this simple example, the solution for $\bM_{\smbbeta}$ is 
$$\bM_{\smbbeta}=
\frac{\phi}{a_1\big(\betazZero,\betaoZero,\sigsqZero\big)}
\left[
\begin{array}{cc}
a_2\big(\betaoZero,\sigsqZero\big)   &\quad -E\big(Xe^{\betaoZero X}\big) \\[1ex]
-E\big(Xe^{\betaoZero X}\big)        &\quad E\big(e^{\betaoZero X}\big)
\end{array}
\right].
$$
Studentisation of the two-term asymptotic covariance matrix for obtaining confidence
intervals and Wald hypothesis tests is straightforward. For example, $E(X^2 e^{\betaoZero X})$
can be replaced by the estimator
$$\frac{1}{mn}\sumim\sum_{j=1}^{n_i} X_{ij}^2e^{\betahat_1 X_{ij}}.$$
This practical aspect is discussed in depth in Section \ref{sec:utility}.

The remainder of this section is concerned with the \emph{theoretical} problem
of obtaining the forms of $\bM_{\smbbeta}$ and $\bM_{\smbSigma}$ for 
model (\ref{eq:theModel}) in general. The achievement of this goal has turned
out to be quite challenging. The score asymptotic approximation approach
used in Jiang \textit{et al.} (2022) requires higher numbers of terms to obtain
valid two-term covariance matrix approximations. Some of these terms can only be expressed 
using three-dimensional arrays rather than with matrices. Succinct statement of
$\bM_{\smbbeta}$ and $\bM_{\smbSigma}$ is only possible with well-designed
nested function notation. A novel notation for multiplicative combining of 
three-dimensional arrays with compatible matrices is also beneficial. 
The next subsection focusses on these notational aspects.

\subsection{Notation for the Main Result}\label{sec:theNotation}

Let $\bAsc$ be a $d_1\times d_2\times d_3$ array and $\bM$ be a $d_1\times d_2$
matrix. Then we let
\begin{equation}
\bAsc\bigstar\bM\quad\mbox{denote the $d_3\times1$ vector with $t$th entry given by}
\quad\sum_{r=1}^{d_1}\sum_{s=1}^{d_2}(\bAsc)_{rst}(\bM)_{rs}.
\label{eq:BlackStarDefn}
\end{equation}
Next, for $\bU\sim N(\bzero,\bSigmaZero)$, define
$$\bOmegaAA(\bU)\equiv E\Big\{b''\big((\bbetaAzero+\bU)^T\bXA+(\bbetaBzero)^T\bXB\big)\bXA\bXA^T|\bU\Big\},$$
$$\bOmegaAB(\bU)\equiv E\Big\{b''\big((\bbetaAzero+\bU)^T\bXA+(\bbetaBzero)^T\bXB\big)\bXA\bXB^T|\bU\Big\}$$
and
$$\bOmegaBB(\bU)\equiv E\Big\{b''\big((\bbetaAzero+\bU)^T\bXA+(\bbetaBzero)^T\bXB\big)\bXB\bXB^T|\bU\Big\}.$$
Also let $\bOmegadAAA(\bU)$ be the $\dR\times\dR\times\dR$ array with $(r,s,t)$ entry equal to 
$$E\Big\{b'''\big(\linXABUzero\big)(\bXA)_r(\bXA)_s(\bXA)_t|\bU\Big\}.$$
and $\bOmegadAAB(\bU)$ be the $\dR\times\dR\times(\dF-\dR)$ array with $(r,s,t)$ entry equal to 
$$E\Big\{b'''\big(\linXABUzero\big)(\bXA)_r(\bXA)_s(\bXB)_t\big|\bU\Big\}.$$
Define the random vectors:
{\setlength\arraycolsep{1pt}
\begin{eqnarray*}
\VECTORone&\equiv&\vech(\bSigma-\bU\bU^T),\quad
\VECTORtwo\equiv \bOmegadAAA(\bU)\bigstar\bOmegaAA(\bU)^{-1},\\[1ex]
\VECTORthree&\equiv&\bOmegadAAB(\bU)\bigstar\bOmegaAA(\bU)^{-1}
\ \ \mbox{and}\ \  
\VECTORfour\equiv\bD_{\dR}^+\vecof\Big(\bOmegaAA(\bU)^{-1}\bSigma^{-1}\big\{\bSigma- \bU\bU^T
-\bSigma\VECTORtwo\bU^T\big\}\Big).
\end{eqnarray*}
}
Then define the random matrices:
{\setlength\arraycolsep{1pt}
\begin{eqnarray*}
\MATRIXfive&\equiv&\bOmegaAA(\bU)^{-1}\bOmegaAB(\bU),\quad 
\MATRIXsix\equiv\bOmegaBB(\bU)-\MATRIXfive^T\bOmegaAB(\bU),\\[1ex]
\MATRIXseven&\equiv&\bU\bU^T\bSigma^{-1}\bOmegaAA(\bU)^{-1},\quad
\MATRIXeight\equiv\bD_{\dR}^+\big[(\bU\bU^T)\otimes\{\bOmegaAA(\bU)^{-1}\}\big]\bD_{\dR}^{+T},\\[1ex]
\quad\mbox{and}\quad
\MATRIXnine&\equiv&\VECTORone\VECTORfour^T+\VECTORfour\VECTORone^T.
\end{eqnarray*}
}
Lastly, define the expectation matrices:
{\setlength\arraycolsep{1pt}
\begin{eqnarray*}
\bLambdaAA&\equiv&E\Big\{\MATRIXseven+\MATRIXseven^T-\bOmegaAA(\bU)^{-1}
+\bOmegaAA(\bU)^{-1}\VECTORtwo\bU^T+\bU\VECTORtwo^T\bOmegaAA(\bU)^{-1}\Big\},\\[1ex]
\bLambdaAB&\equiv&E\Big\{\bU\bU^T\bSigma^{-1}\MATRIXfive+\bU\VECTORtwo^T\MATRIXfive-\bU\VECTORthree^T\Big\}
\quad\mbox{and}\\[1ex]
\bPhi&\equiv&E\Big(\Big[\MATRIXfive^T\big\{\bSigma^{-1}\bU+\VECTORtwo\big\}-\VECTORthree\Big]\VECTORone^T\Big).
\end{eqnarray*}
}

\subsection{Assumptions for the Main Result}\label{sec:theAssums}

The main result depends on the following sample size asymptotic assumptions:
\begin{itemize}
\item[] The number of groups $m$ diverges to $\infty$.
\item[] The within-group sample sizes $n_i$ diverge to $\infty$ in such a way that
$n_i/n\to C_i$ for constants $0<C_i<\infty$, $1\le i\le m$.
\item[] The ratio $n/m$ converges to zero.
\end{itemize}
The last of these conditions is in keeping with the number of groups being
large compared with the within-group sample sizes, as often arises in practice.
For our asymptotics it ensures that, for the harder-to-estimate parameters, 
the asymptotic variances of the maximum likelihood estimators have 
leading terms of the form $C_1m^{-1}+C_2(mn)^{-1}$. In addition, it ensures
that the Fisher information is sufficiently dominant for obtaining
asymptotic variances.

We also assume that the $(\bX,\bU)$ joint distribution is such that
all required convergence in probability limits that appear in the 
deterministic order $(mn)^{-1}$ terms are justified.  
An example of such a convergence in probability statement is 
\begin{equation}
\begin{array}{c}
{\displaystyle\oon}
E\Big(\HscrABi^T\HscrAAi^{-1}\HscrABi\Big|\bX_i\Big)\convprobBka\bOmegaAB(\bU)^T
\bOmegaAA(\bU)^{-1}\bOmegaAB(\bU)\\[1ex]
\mbox{where}\ \HscrAAi\equiv\dispsumjni b''\big(\linXABUi\big)\bXAij\bXAij^T\\[1ex]
\mbox{and}\ \HscrABi\equiv\dispsumjni b''\big(\linXABUi\big)\bXAij\bXBij^T.
\end{array}
\label{eq:NormaTriangle}
\end{equation}
Assumption (A3) of Jiang \textit{et al.} (2022) provides a
moment-type condition that is sufficient for (\ref{eq:NormaTriangle}) 
to hold. Also, we assume that the tail behaviour of the $(\bX,\bU)$
distribution is such that statements concerning the $\myoP\{(mn)^{-1}\}$
remainder terms are valid. The determination of sufficient conditions on
the $(\bX,\bU)$ distribution that guarantee the validity of the main
result is a tall order, and beyond the scope of this article.

\subsection{Statement of the Main Result}\label{sec:actualMainResult}

Using the notation presented in Section \ref{sec:theNotation} and
under the assumptions described in Section \ref{sec:theAssums},
and assuming $\dF>\dR$ we have
\begin{equation}
{\setlength\arraycolsep{0pt}
\begin{array}{l}
\Cov\big(\bbetahat|\bXsc\big)=
\displaystyle{\frac{1}{m}}
\left[
\begin{array}{cc}
\ \bSigmaZero &\ \ \ \bO \\[1ex]
\bO     &\ \ \ \bO
\end{array}
\right]
+
\displaystyle{\frac{\phi}{mn}}
\left[
\begin{array}{cc}
\bLambdaAA^{-1}             &\quad\bLambdaAA^{-1}\bLambdaAB \\[2ex]
\bLambdaAB^T\bLambdaAA^{-1} &\quad\bLambdaAB^T\bLambdaAA^{-1}\bLambdaAB
                              +E\big\{\MATRIXsix\big\}
\end{array}
\right]^{-1}\{1+\myoP(1)\}\quad\mbox{and}\\[6ex]
\Cov\big(\vech(\bSigmaMLE)|\bXsc\big)=
\displaystyle{\frac{2\bD_{\dR}^+(\bSigmaZero\otimes\bSigmaZero)
\bD_{\dR}^{+T}}{m}}\\[2ex]
\qquad\qquad\qquad\qquad\quad\qquad
+\displaystyle{\frac{\phi}{mn}}\tcbk{\Big(}2E\big\{\MATRIXnine-2\MATRIXeight\big\}
+\bPhi^T\tcbk{\big[}E\{\MATRIXsix\}\tcbk{\big]}^{-1}\bPhi\tcbk{\Big)}\{1+\myoP(1)\}.
\end{array}
}
\label{eq:mainResult}
\end{equation}
For the $\dF=\dR$ boundary case the first term of $\Cov\big(\bbetahat|\bXsc\big)$
is simply $\frac{1}{m}\bSigmaZero$.

A supplement to this article contains a full derivation of (\ref{eq:mainResult}).

\subsubsection{The Gaussian Response Special Case}

In the Gaussian response special case we have $b''(x)=1$ and
$b'''(x)=0$ and the main result reduces to the following
succinct form:
\begin{equation}
{\setlength\arraycolsep{1pt}
\begin{array}{rcl}
\Cov\big(\bbetaMLE|\bXsc\big)&=&
\displaystyle{\frac{1}{m}}
\left[
\begin{array}{cc}
\bSigmaZero &\quad \bO \\[1ex]
\bO     &\quad  \bO
\end{array}
\right]
+\displaystyle{\frac{\phi\big\{E(\bX\bX^T)\big\}^{-1}\{1+\myoP(1)\}}{mn}}\quad\mbox{and}\\[5ex]
\Cov\big(\vech(\bSigmaMLE)|\bXsc\big)
&=&\displaystyle{\frac{2\bD_{\dR}^+(\bSigmaZero\otimes\bSigmaZero)\bD_{\dR}^{+T}}{m}}\\[2ex]
&&\qquad\qquad
+\displaystyle{\frac{4\phi\bD_{\dR}^+\big[\bSigmaZero\otimes\{E(\bXA\bXA^T)\}^{-1}\big]\bD_{\dR}^{+T}
\{1+\myoP(1)\}}{mn}}.
\end{array}
}
\label{eq:mainResultGauss}
\end{equation}

\noindent
We are not aware of any previous appearances of (\ref{eq:mainResultGauss}) in the wider 
linear mixed model literature.

\section{Utility of the Second Term Improvements}\label{sec:utility}

We now describe the utility of (\ref{eq:mainResult}) in statistical
contexts such as inference and design. Improved confidence intervals
is \tcbk{a} particularly straightforward application, which we treat next.

\subsection{Confidence Intervals}\label{sec:CIconstruct}

For any $\bu\in\real^{\dR}$, define
\begin{equation}
{\setlength\arraycolsep{1pt}
\begin{array}{rcl}
\bOmegaAAhat(\bu)&\equiv& 
{\displaystyle\frac{1}{mn}\sumim\sumjni}
b''\Big((\bbetaAMLE+\bu)^T\bXAij+\bbetaBMLE^T\bXBij\Big)\bXAij\bXAij^T,\\[1ex]
\bOmegaABhat(\bu)&\equiv& 
{\displaystyle\frac{1}{mn}\sumim\sumjni}
b''\Big((\bbetaAMLE+\bu)^T\bXAij+\bbetaBMLE^T\bXBij\Big)\bXAij\bXBij^T\\[1ex]
\mbox{and}\quad
\bOmegaBBhat(\bu)&\equiv& 
{\displaystyle\frac{1}{mn}\sumim\sumjni}
b''\Big((\bbetaAMLE+\bu)^T\bXAij+\bbetaBMLE^T\bXBij\Big)\bXBij\bXBij^T.
\end{array}
}
\label{eq:OmegaHatDefns}
\end{equation}
Then the natural studentisation of $E\big\{\MATRIXsix\big\}$ is
\begin{equation}
{\setlength\arraycolsep{1pt}
\begin{array}{rcl}
\EpsiSixHat&\equiv& 
E\Big\{\bOmegaBBhat(\bU)-\bOmegaABhat(\bU)^T\bOmegaAAhat(\bU)^{-1}
\bOmegaABhat(\bU)\big|\bXsc\Big\}\\[2ex]
&=&|2\pi\bSigmaMLE|^{-1/2}{\displaystyle\int_{\dR}}
\Big\{\bOmegaBBhat(\bu)-\bOmegaABhat(\bu)^T\bOmegaAAhat(\bu)^{-1}
\bOmegaABhat(\bu)\Big\}\exp\Big(-\smhalf\bu^T\bSigmaMLE\bu\Big)\,d\bu.
\end{array}
}
\label{eq:studMain}
\end{equation}
In the last expression of (\ref{eq:studMain}) integration is applied element-wise to each 
entry of the matrix inside the integral.
The natural studentisations of 
\begin{equation}
\bLambdaAA,\quad\bLambdaAB,\quad\bPhi,\quad E\big\{\MATRIXeight\big\}
\quad\mbox{and}\quad E\big\{\MATRIXnine\big\}
\label{eq:TuringMovie}
\end{equation}
are analogous to that for $E\big\{\MATRIXsix\big\}$. 
The studentisations for the quantities in (\ref{eq:TuringMovie}) depend on the functions defined by 
(\ref{eq:OmegaHatDefns}) as well as similar sample counterparts of $\bOmegadAAA(\bU)$ 
and $\bOmegadAAB(\bU)$.
Next define
\begin{equation}
\AsyCovHat(\bbetaMLE)
=\displaystyle{\frac{1}{m}}
\left[
\begin{array}{cc}
\ \bSigmaMLE &\ \ \ \bO \\[1ex]
\bO     &\ \ \ \bO
\end{array}
\right]
+
\displaystyle{\frac{\phihat}{mn}}
\left[
\begin{array}{cc}
\bLambdaAAhat^{-1}             &\quad\bLambdaAAhat^{-1}\bLambdaABhat \\[2ex]
\bLambdaABhat^T\bLambdaAAhat^{-1} &\quad\bLambdaABhat^T\bLambdaAAhat^{-1}\bLambdaABhat
                              +\EpsiSixHat
\end{array}
\right]^{-1}
\label{eq:WarrenTheWizard}
\end{equation}
and
\begin{equation}
\begin{array}{l}
\AsyCovHat\big(\vech(\bSigmaMLE)\big)
=\displaystyle{\frac{2\bD_{\dR}^+(\bSigmahat\otimes\bSigmahat)
\bD_{\dR}^{+T}}{m}}\\[1ex]
\quad\ \qquad\qquad\qquad\qquad
+\displaystyle{\frac{\phihat}{mn}}\Big(2\EpsiNineHat
-4\EpsiEightHat
+\bPhihat^T\big[\EpsiSixHat\big]^{-1}\bPhihat\Big).
\end{array}
\label{eq:PaigeTheBabe}
\end{equation}
In the general quasi-likelihood situation, the most common choice for $\phihat$ 
is the method of moments estimator and is often labelled the \emph{Pearson} estimator. 
For ordinary likelihood settings, such as for Gaussian and Gamma responses, 
$\phihat$ could instead be the maximum likelihood estimator.

Let $(\bbetaZero)_k$ denote the $k$th entry of $\bbetaZero$. Then approximate 
$100(1-\alpha)\%$ confidence intervals for $(\bbetaZero)_k$ based
on (\ref{eq:WarrenTheWizard}) are 
\begin{equation}
(\bbetaMLE)_k\pm\Phi^{-1}(1-\smhalf\alpha)\sqrt{\left\{\AsyCovHat(\bbetaMLE)\right\}_{kk}},
\quad 1\le k\le \dF.
\label{eq:betaCI}
\end{equation}
The confidence intervals in (\ref{eq:betaCI}) are analogous to those given
in Section 4 of Jiang \textit{et al.} (2022). For $1\le k\le\dR$, 
(\ref{eq:betaCI}) provides second term improvements of the Jiang \textit{et al.} (2022)
confidence intervals. For $\dR+1\le k\le \dF$ both sets of confidence intervals
are identical. 

Improved confidence intervals for the random effects
covariance parameters can be constructed in a similar fashion
based on (\ref{eq:PaigeTheBabe}).

\subsection{Other Utilities}

The second term improvements of (\ref{eq:mainResult}) may also be 
applied to Wald hypothesis tests and sample size calculations. 
Optimal design is another possible utility, but would require
second term improvements of the type of theory given in 
Section 5 of Jiang \textit{et al.} (2022).

\section{Numerical Results}\label{sec:numericRes}

We conducted a simulation exercise aimed at understanding potential 
practical impacts of second term improvements to generalized linear
mixed model asymptotics. The results are presented in this section.

Our simulation exercise involved generation of data sets from the 
$\dF=5$ and $\dR=2$ logistic mixed model
\begin{equation}
\begin{array}{l}
Y_{ij}|X_{1ij},X_{2ij},X_{3ij},X_{4ij},U_i\ \mbox{independently distributed as}\\[1ex]
\qquad
\mbox{Bernoulli}
\Big(1/\tcbk{\big(1+}\exp\tcbk{[}-\{\betaNoughtZero+U_{0i}+(\betaOneZero+U_{1i})X_{1ij}\tcbk{+}
\betaTwoZero X_{2ij}+\betaThreeZero X_{3ij}+\betaFourZero X_{4ij}
\}\tcbk{]\big)}\Big),\\[1ex]
\mbox{where the}\ 
\left[
\begin{array}{c}
U_{0i}\\
U_{1i}
\end{array}
\right]\ \mbox{are independent}\ N(\bzero,\bSigmaZero)\ \mbox{random vectors},
\ 1\le i\le m,\ 1\le j\le n.
\end{array}
\label{eq:bigBird}
\end{equation}
The `true' parameter values were set to 
\begin{equation}
\big(\betaNoughtZero,\betaOneZero,\betaTwoZero,\betaThreeZero,\betaFourZero\big)
=(0.35,0.96,-0.47,1.06,-1.31)\quad\mbox{and}\quad 
\bSigmaZero=\left[
\begin{array}{cc}
0.56  & -0.34\\[1ex]
-0.34 & 0.89
\end{array}
\right]
\label{eq:littleWing}
\end{equation}
and the predictor data were generated from independent Uniform distributions
on the unit interval. To assess potential large sample improvements afforded
by the two-term asymptotic covariance expressions at (\ref{eq:mainResult})
we varied $m$ over the set $\{100,150,\ldots,500\}$ and fixed $n$ at $m/10$.
For each $(m,n)$ pair we then simulated $500$ data sets
according to (\ref{eq:bigBird}) and (\ref{eq:littleWing}) and obtained 
approximate 95\% confidence intervals for all model parameters according
to the approach described in Section 4 of Jiang \textit{et al.} (2022)
and the second term improvements described in Section 
\ref{sec:CIconstruct} of this article. The requisite bivariate
integrals were obtained using the function \texttt{hcubature()} within
the \textsf{R} language package \textsf{cubature} 
(Balasubramanian \textit{et al.}, 2023).

\begin{figure}[!ht]
\includegraphics[width=\textwidth]{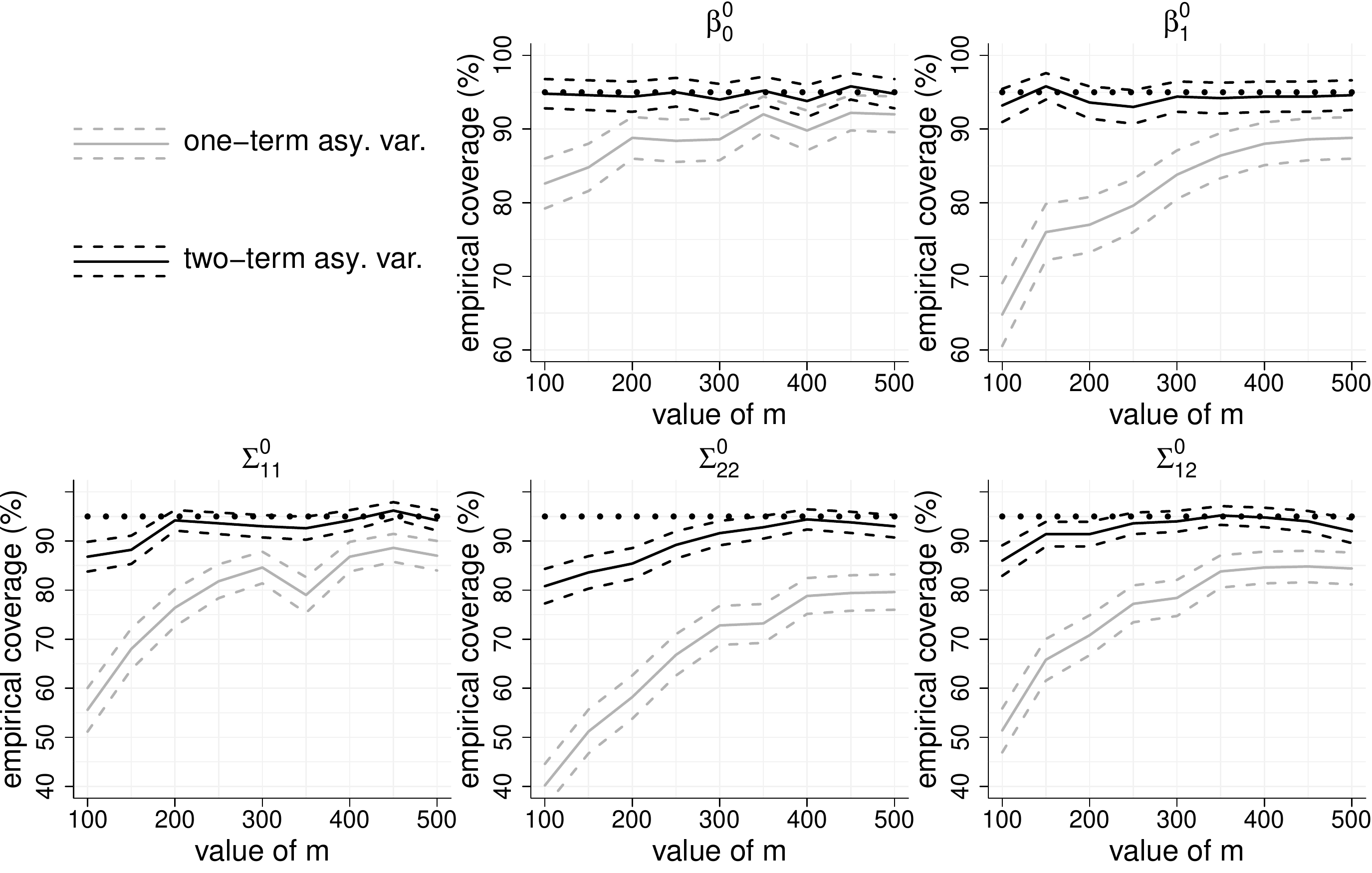}
\caption{\it Empirical coverage of confidence intervals from the simulation
exercise described in the text. Each panel corresponds to a model
parameter that is impacted by second term asymptotic improvements.
The advertised coverage level is fixed at 95\% and is indicated by \tcbk{a}
horizontal dotted line in each panel. The solid curves show, dependent on 
the number of groups $m$, the empirical coverage levels for confidence 
intervals that use both one-term and two-term asymptotic variance approximations. 
The dashed curves correspond to plus and minus two standard errors of
the sample proportions. The within-group sample size, $n$, is fixed at $m/10$.
}
\label{fig:MBWcovMainFig} 
\end{figure}

Note that the confidence intervals for $\betaNoughtZero$, $\betaOneZero$ 
and the entries of $\bSigmaZero$ differ according to the two approaches since the estimators
of these parameters have order $m^{-1}$ asymptotic variances. The confidence
intervals for $\betaTwoZero$, $\betaThreeZero$ and $\betaFourZero$
are unaffected by the second term asymptotic improvements since
their estimators have order $(mn)^{-1}$ asymptotic variances.

Figure \ref{fig:MBWcovMainFig} compares the empirical coverages of confidence intervals
with advertised levels of 95\% for the one-term asymptotic variances of 
Jiang \textit{et al.} (2022) and the two-term asymptotic variances that
arise from (\ref{eq:mainResult}). In Figure \ref{fig:MBWcovMainFig} we
only consider the parameters that are affected by second term improvement.
\tcbk{The empirical coverages for the other parameters are provided in the supplement.}

It is clear from Figure \ref{fig:MBWcovMainFig} that our second term improvements
lead to much better coverages for lower sample size situations. On the other
hand, one-term confidence intervals are trivial to compute whilst the two-term
versions require considerable computing involving numerical integration. 

Simulation results such as those summarised by Figure \ref{fig:MBWcovMainFig}
provide an appreciation for the practical trade-offs 
arising from precise asymptotics for \mygeneralized\ linear mixed models.

\section*{Acknowledgements}

\tcbk{We are grateful to Alessandra Salvan and Nicola Sartori for advice related to this research.} 
This research was supported by the Australian Research Council Discovery
Project DP230101179.

\section*{References}

\bib
Balasubramanian, N., Johnson, S.G., Hahn, T., Bouvier, A. \myand Ki{\^e}u, K. (2023).
\textsf{cubature 2.0.4.6}: Adaptive multivariate integration over hypercubes.
\textsf{R} package.\\ 
\texttt{https://r-project.org}

\bib
Bhaskaran, A. and Wand, M.P. (2023).
Dispersion parameter extension of precise
generalized linear mixed model asymptotics.
\textit{Statistics and Probability Letters}, 
\textbf{193}, Article 109691.

\bib
Jiang, J. \myand Nguyen, T. (2021). \textit{Linear and Generalized Linear 
Mixed Models and Their Applications, Second Edition.} New York: Springer.

\bib
Jiang, J., Wand, M.P. \myand Bhaskaran, A. (2022).
Usable and precise asymptotics for generalized linear mixed model analysis
and design. \textit{Journal of the Royal Statistical Society, Series B},
\textbf{84}, 55--82.

\bib
McCulloch, C.E., Searle, S.R. \myand Neuhaus, J.M. (2008). 
\textit{Generalized, Linear, and Mixed Models. Second Edition.}
New York: John Wiley \& Sons.

\bib
Stroup, W.W. (2013). \textit{Generalized Linear Mixed Models.}
Boca Raton, Florida: CRC Press.

\vfill\eject
%
%
\renewcommand{\theequation}{S.\arabic{equation}}
\renewcommand{\thesection}{S.\arabic{section}}
\renewcommand{\thetable}{S.\arabic{table}}
\renewcommand{\thefigure}{S.\arabic{figure}}
\setcounter{equation}{0}
\setcounter{table}{0}
\setcounter{section}{0}
\setcounter{figure}{0}
\setcounter{page}{1}
\setcounter{footnote}{0}

\begin{center}

{\Large Supplement for:}
\vskip3mm

\centerline{\Large\bf Second Term Improvements to Generalised}
\vskip1mm
\centerline{\Large\bf Linear Mixed Model Asymptotics}
\vskip7mm
\centerline{\normalsize\sc By Luca Maestrini$\null^1$, 
Aishwarya Bhaskaran$\null^2$ and Matt P. Wand$\null^3$}
\vskip5mm
\centerline{\textit{$\null^1$\tcbk{The} Australian National University, 
$\null^2$Macquarie University 
and 
$\null^3$University of Technology Sydney}}
\vskip6mm
\end{center}

\section{Introduction}

The purpose of this supplement is to provide detailed derivational steps for
the main result of Section \ref{sec:actualMainResult} \tcbk{and further
details on our simulation exercise}. 
Sections \ref{sec:matAlg}--\ref{sec:Miyata} provide
relevant results concerning matrix algebra and multivariate calculus.
In Sections \ref{sec:exactScores}--\ref{sec:scoreApprox} we focus on the
scores of the model parameters and their high-order asymptotic approximations.
Sections \ref{sec:outProd} and \ref{sec:FishInfo} are concerned with approximation 
of the Fisher information matrix. The final stages of the derivations of
(\ref{eq:mainResult}) and (\ref{eq:mainResultGauss}) are given in
Sections \ref{sec:ApproxAsyCov} and \ref{sec:GaussCase}. \tcbk{Section \ref{sec:additFig}
provides some additional empirical coverage plots from the logistic mixed model simulation
exercise described in Section \ref{sec:numericRes}.}

\section{Matrix Algebraic Results}\label{sec:matAlg}

The derivation of the results in Section \ref{sec:actualMainResult} benefits from 
particular matrix results, which are summarized in this section.

For each $d\in\naturalNumbers$ the $d^2\times\smhalf d(d+1)$ matrix $\bD_d$ and 
$d^2\times d^2$ matrix $\bK_d$ are constant matrices containing zeroes and ones 
such that 
$$\bD_d\vech(\bA)=\vecof(\bA)\quad\mbox{for all symmetric $d\times d$ matrices $\bA$} $$
and
$$\bK_d\vecof(\bB)=\vecof(\bB^T)\quad\mbox{for all $d\times d$ matrices $\bB$}.$$
Examples are 
$$\bD_2=\left[
\begin{array}{ccc}
1 & 0 & 0\\
0 & 1 & 0\\
0 & 1 & 0\\
0 & 0 & 1
\end{array}
\right]\quad\mbox{and}\quad
\bK_2=\left[
\begin{array}{cccc}
1 & 0 & 0 & 0\\
0 & 0 & 1 & 0\\
0 & 1 & 0 & 0\\
0 & 0 & 0 & 1
\end{array}
\right].
$$
The $\bD_d$ are called \emph{duplication} matrices, whilst the $\bK_d$ are
called \emph{commutation} matrices. As stated in Section \ref{sec:modelDescrip}, the 
Moore-Penrose inverse of $\bD_d$ is $\bD_d^+=(\bD_d^T\bD_d)^{-1}\bD_d^T$.
Chapter 3 of Magnus \myand Neudecker (1999)
contains several results concerning these families of matrices, a few of which
are relevant to the derivation of (\ref{eq:mainResult}). For convenience, we list 
them here. 

Theorem 9(c) in Chapter 3 of Magnus \myand Neudecker (1999) implies that for
any $d\times d$ matrix $\bA$ and $d\times 1$ vector $\bb$, we have
\begin{equation}
\bK_d(\bA\otimes\bb)=\bb\otimes\bA.
\label{eq:MandNone}
\end{equation}
Theorem 12(a) in the same chapter asserts that 
\begin{equation}
\bK_d\bD_d=\bD_d
\label{eq:MandNtwo}
\end{equation}
and implies that, for any $d\times d$ matrix $\bA$,
\begin{equation}
\bD_d^T\vecof(\bA)=\bD_d^T\vecof(\bA^T).
\label{eq:MandNthree}
\end{equation}
Also, Theorem 13(b) and Theorem 13(d) provide for a $d\times d$ matrix $\bA$
\begin{equation}
\bD_d\bD_d^+(\bA\otimes\bA)\bD_d^{+T}=(\bA\otimes\bA)\bD_d^{+T}
\label{eq:MandNfour}
\end{equation}
and, assuming that $\bA$ is invertible,
\begin{equation}
\big\{\bD_d^T(\bA\otimes\bA)\bD_d\big\}^{-1}=\bD_d^+(\bA^{-1}\otimes\bA^{-1})\bD_d^{+T}.
\label{eq:MandNfive}
\end{equation}
Lastly, we state \tcbk{two} matrix identit\tcbk{ies} that are used in
the derivations. For matrices $\bA$, 
$\bB$ and $\bC$ such that $\bA\bB\bC$ is defined, we have
\begin{equation}
\vecof(\bA\bB\bC)=(\bC^T\otimes\bA)\vecof(\bB).
\label{eq:vecABC}
\end{equation}
\tcbk{For conformable matrices $\bA$, $\bB$, $\bC$ and $\bD$, we have}
\begin{equation}
\tcbk{(\bA\otimes\bB)(\bC\otimes\bD)=(\bA\bC)\otimes(\bB\bD).
\label{eq:ACkronBD}}
\end{equation}

\section{Multivariate Derivative Notation}

For $f$ a smooth real-valued function of the $d$-variate argument 
$\bx\equiv(x_1,\ldots,x_d)$, let $\nabla f(\bx)$ denote
the $d\times1$ vector with $r$th entry $\partial f(\bx)/\partial x_r$,
$\nabla^2 f(\bx)$ denote the $d\times d$ matrix with $(r,s)$ entry 
$\partial^2 f(\bx)/(\partial x_r\partial x_s)$ and 
$\nabla^3 f(\bx)$ denote the $d\times d\times d$ array with $(r,s,t)$ entry 
$\partial^3 f(\bx)/(\partial x_r\partial x_s\partial x_t)$.

\section{Three-Term Taylor Series Expansion of Gradient Vectors}

Consider $f:\real^d\to\real$. Then, for sufficiently smooth $f$, the three-term Taylor 
expansion of 
$$f(\bx+\bh)\quad\mbox{where $\bx\equiv(x_1,\ldots,x_d)$ and $\bh\equiv(h_1,\ldots,h_d)$}$$
is
\begin{equation}
f(\bx+\bh)=f(\bx)+\sum_{r=1}^d\{\nabla f(\bx)\}_rh_r
+\smhalf\sum_{r=1}^d\sum_{s=1}^d\{\nabla^2 f(\bx)\}_{rs}h_rh_s+\ldots\tcbk{.}
\label{eq:MenloPark}
\end{equation}
Now consider $\alpha:\real^d\to\real$ and its gradient function $\nabla\alpha:\real^d\to\real^d$.
If (\ref{eq:MenloPark}) is applied to each entry of $(\nabla\alpha)(\bx+\bh)$ then we have
$$
\left[
{\setlength\arraycolsep{0.5pt}
\begin{array}{c}
\{(\nabla\alpha)(\bx+\bh)\}_1\\
\vdots\\
\{(\nabla\alpha)(\bx+\bh)\}_d
\end{array}
}
\right]
=\left[{\setlength\arraycolsep{0.5pt}
\begin{array}{c}
\{(\nabla\alpha)(\bx)\}_1+{\displaystyle\sum_{r=1}^d}\{\nabla^2\alpha(\bx)\}_{r1}h_r
+\smhalf{\displaystyle\sum_{r=1}^d\sum_{s=1}^d}\{\nabla^3\alpha(\bx)\}_{rs1}h_rh_s\\
\vdots\\
\{(\nabla\alpha)(\bx)\}_{\tcbk{d}}+{\displaystyle\sum_{r=1}^d}\{\nabla^2\alpha(\bx)\}_{r\tcbk{d}}h_r
+\smhalf{\displaystyle\sum_{r=1}^d\sum_{s=1}^d}\{\nabla^3\alpha(\bx)\}_{rsd}h_rh_s\
\end{array}
}
\right]+\ldots.
$$
From this it is clear that
\begin{equation}
(\nabla\alpha)(\bx+\bh)=(\nabla\alpha)(\bx)+
\{(\nabla^2\alpha)(\bx)\}\bh+\smhalf
\{(\nabla^3\alpha)(\bx)\}\bigstar(\bh\bh^T)+\ldots.
\label{eq:TuskVideo}
\end{equation}
where the $\bigstar$ notation is as defined by 
(\ref{eq:BlackStarDefn}).

\section{Higher Order Approximation of Multivariate Integral Ratios}\label{sec:Miyata}

The main tool for approximation of the Fisher information matrix
of (\ref{eq:theModel}) is higher order Laplace-type approximation of 
multivariate integral ratios. Appendix A of Miyata (2004) provides
such a result, which states that for smooth real-valued $d$-variate 
functions $g$, $c$ and $h$\tcbk{,}
\begin{equation}
{\setlength\arraycolsep{1pt}
\begin{array}{rcl}
&&{\displaystyle\frac{\int_{\real^d}g(\bx)c(\bx)\exp\{-n h(\bx)\}\,d\bx}
{\int_{\real^d} c(\bx)\exp\{-n h(\bx)\}\,d\bx}}
=g(\bx^*)+
{\displaystyle\frac{\nabla g(\bx^*)^T\{\nabla^2 h(\bx^*)\}^{-1}\nabla c(\bx^*)}{n c(\bx^*)}}\\[3ex]
&&\qquad+{\displaystyle\frac{\tr[\{\nabla^2 h(\bx^*)\}^{-1}\nabla^2 g(\bx^*)]}{2n}}
-{\displaystyle\frac{\nabla g(\bx^*)^T\{\nabla^2 h(\bx^*)\}^{-1}
\Big[\nabla^3 h(\bx^*)\bigstar\{\nabla^2 h(\bx^*)\}^{-1}\Big]}
{2n}}+O(n^{-2})
\end{array}
}
\label{eq:Miyata}
\end{equation}
where 
$$\bx^*\equiv\argmin{\bx\in\real^d}h(\bx).$$

\section{Exact Score Expressions}\label{sec:exactScores}

For $1\le i\le m$, let $\pYiGXi$ denote the conditional density function, 
or probability mass function, of $\bY_i$ given $\bX_i$. Then let
$$\bSAi\equiv\nablabbetaA\log\pYiGXi(\bY_i|\bX_i),
\quad \bSBi\equiv\nablabbetaB\log\pYiGXi(\bY_i|\bX_i)
$$
and
$$\bSCi\equiv\nablavechbSigma\log\pYiGXi(\bY_i|\bX_i)$$
denote the $i$th contribution to the scores with respect to 
each of $\bbetaA$, $\bbetaB$ and $\vech(\bSigma)$. 
Then it is straightforward to show that the exact scores are
\begin{equation}
\bSAi=\frac{\intdR\bgiA(\bu)\cS(\bu)\exp\{-n h_i(\bu)\}\,d\bu}{\intdR \cS(\bu)\exp\{-n h_i(\bu)\}\,d\bu},
\label{eq:scoreA}
\end{equation}
\begin{equation}
\bSBi=\frac{\intdR\bgiB(\bu)\cS(\bu)\exp\{-n h_i(\bu)\}\,d\bu}{\intdR \cS(\bu)\exp\{-n h_i(\bu)\}\,d\bu}
\label{eq:scoreB}
\end{equation}
and
\begin{equation}
\bSCi=\frac{\intdR \bgiC(\bu)\cS(\bu)\exp\{-n h_i(\bu)\}\,d\bu}{\intdR \cS(\bu)\exp\{-n h_i(\bu)\}\,d\bu}
-\smhalf\bD_{\dR}^T\vecof(\bSigma^{-1})
\label{eq:scoreC}
\end{equation}
where 
{\setlength\arraycolsep{1pt}
\begin{eqnarray*}
\cS(\bu)&\equiv&\exp(-\smhalf\bu^T\bSigma^{-1}\bu),\quad \bgiA(\bu)\equiv\bSigma^{-1}\bu,\\[1ex]
\bgiB(\bu)&\equiv&\frac{1}{\phi}\dispsumjn\bXBij\{Y_{ij}-b'\big(\linXABu\big)\},\\[1ex]
\bgiC(\bu)&\equiv&\smhalf\bD_{\dR}^T(\bSigma^{-1}\otimes\bSigma^{-1})\vecof(\bu\bu^T)\\[1ex]
\mbox{and}\quad
h_i(\bu)&\equiv&\,-\frac{1}{n\phi}\dispsumjn\left\{Y_{ij}\bu^T\bXAij
-b\big(\linXABu\big)\right\}.
\end{eqnarray*}
}
An integration by parts step is used to obtain the $\bSAi$ expression.

In the upcoming sections we obtain asymptotic approximations of $\bSAi$
$\bSBi$ and $\bSCi$. Key quantities for these approximations are
$$\bU_i^*\equiv\argmin{\bu\in\real^{\dR}}h_i(\bu),\quad 1\le i\le m.$$

\section{Definitions of Key Summation Quantities}

Our derivation of (\ref{eq:mainResult}) involves manipulations of particular
summation quantities, which are defined in this section. At the 
end of this section we state some important moment-type relationships between
the quantities.

For each $1\le i\le m$, define $\GscrAi$, $\GscrBi$, $\HscrAAi$, $\HscrABi$, 
and $\HscrBBi$ as follows:
{\setlength\arraycolsep{1pt}
\begin{eqnarray*}
\GscrAi&\equiv&\dispsumjni\{Y_{ij} -b'\big(\linXABUi\big)\}\bXAij,\\[1ex]
\GscrBi&\equiv&\dispsumjni\{Y_{ij} -b'\big(\linXABUi\big)\}\bXBij,\\[1ex]
\HscrAAi&\equiv&\dispsumjni b''\big(\linXABUi\big)\bXAij\bXAij^T,
\end{eqnarray*}
}
{\setlength\arraycolsep{1pt}
\begin{eqnarray*}
\HscrABi&\equiv&\dispsumjni b''\big(\linXABUi\big)\bXAij\bXBij^T\\[1ex]
\mbox{and}\quad
\HscrBBi&\equiv&\dispsumjni b''\big(\linXABUi\big)\bXBij\bXBij^T.
\end{eqnarray*}
}

In a similar vein, define $\HscrdAAAi$ to be the $\dR\times\dR\times\dR$ 
array with $(r,s,t)$ entry equal to 
$$\dispsumjn b'''\big(\linXABUi\big)(\bXAij)_r(\bXAij)_s(\bXAij)_t$$
and $\HscrdAABi$ to be the $\dR\times\dR\times\dB$ array with $(r,s,t)$ entry 
equal to 
$$\dispsumjni b'''\big(\linXABUi\big)(\bXAij)_r(\bXAij)_s(\bXBij)_t$$
where
$$\dB\equiv \dF-\dR.$$

The following relationships are of fundamental importance for the 
derivation of (\ref{eq:mainResult}):
\begin{equation}
\begin{array}{c}
E\big(\GscrAi|\bX_i,\bU_i\big)=\bzero,\quad E\big(\GscrBi|\bX_i,\bU_i\big)=\bzero,\\[1ex]
E\big(\GscrAi^{\otimes2}|\bX_i,\bU_i\big)=\phi\HscrAAi,\quad
E\big(\GscrAi\GscrBi^T|\bX_i,\bU_i\big)=\phi\HscrABi
\tcbk{\quad\mbox{and}\quad E\big(\GscrBi^{\otimes2}|\bX_i,\bU_i\big)=\phi\HscrBBi}
\end{array}
\label{eq:Ulverstone}
\end{equation}
where, throughout this supplement,
$$\bv^{\otimes2}\equiv\bv\bv^T\quad\mbox{for any column vector $\bv$}.$$
Also note that 
$$\GscrAi=\myOP(n^{1/2})\bone_{\dR},\quad \GscrBi=\myOP(n^{1/2})\bone_{\dB},\quad
\HscrAAi=\myOP(n)\bone_{\dR}^{\otimes2},
$$
$$\HscrBBi=\myOP(n)\bone_{\dB}^{\otimes2},\quad
\HscrABi=\myOP(n)\bone_{\dR}\bone_{\dB}^T
$$
and that all entries of $\HscrdAAAi$ and $\HscrdAABi$ are $\myOP(n)$.

\section{Approximation of $\bU_i^*$}

Use of (\ref{eq:Miyata}) to approximate $\bSAi$, $\bSBi$ and $\bSCi$ 
requires approximation of $\bU_i^*$. Introduce the notation
$\Csc_i(\bu)\equiv n\phi h_i(\bu)$. Then $\bU_i^*$ satisfies
$$\nabla\Csc_i(\bU_i^*)=\bzero$$
where
$$\nabla\Csc_i(\bu)\equiv -\dispsumjni\Big\{Y_{ij}-b'\big(\linXABu\big)\Big\}\bXAij.$$
Then, from (\ref{eq:TuskVideo}) we have
{\setlength\arraycolsep{1pt}
\begin{eqnarray*}
\nabla\Csc_i(\bU_i^*)&=&\nabla\Csc_i(\bU_i+\bU_i^*-\bU_i)\\[1ex]
&=&\nabla\Csc_i(\bU_i)+\{\nabla^2\Csc_i(\bU_i)\}(\bU_i^*-\bU_i)\\
&&\qquad
+\smhalf\left\{\nabla^3\Csc_i(\bU_i)\right\}\bigstar\{(\bU_i^*-\bU_i)(\bU_i^*-\bU_i)^T\}+\ldots.
\end{eqnarray*}
}
Next we seek explicit expressions for $\nabla^2\Csc_i(\bu)$ and $\nabla^3\Csc_i(\bu)$.
Standard vector calculus arguments lead to 
{\setlength\arraycolsep{1pt}
\begin{eqnarray*}
\nabla^2\Csc_i(\bu)&=& \dispsumjni b''\big(\linXABu\big)\bXAij\bXAij^T\\[1ex]
&=&\left[ \dispsumjni b''\big(\linXABu\big)(\bXAij)_r(\bXAij)_s\right]_{1\le r,s\le\dR}.
\end{eqnarray*}
}
Then, the three-dimension array of all third order partial derivatives of $\Csc_i(u)$ is
$$\nabla^3\Csc_i(\bu)=
\left[\dispsumjni b'''\big(\linXABu\big)(\bXAij)_r(\bXAij)_s(\bXAij)_t\right]_{1\le r,s,t\le\dR}.
$$
We then have
$$\nabla\Csc_i(\bU_i^*)
=-\GscrAi+\HscrAAi(\bU_i^*-\bU_i)
+\smhalf\HscrdAAAi\bigstar\{(\bU_i^*-\bU_i)(\bU_i^*-\bU_i)^T\}+\ldots
$$
and so $\nabla\Csc_i(\bU_i^*)=\bzero$ is equivalent to 
\begin{equation}
\HscrAAi^{-1}\GscrAi=(\bU_i^*-\bU_i)+\smhalf\HscrAAi^{-1}
\tcbk{\Big[}\HscrdAAAi\bigstar\{(\bU_i^*-\bU_i)(\bU_i^*-\bU_i)^T\}\tcbk{\Big]}+\ldots.
\label{eq:MariasSister}
\end{equation}
We now invert (\ref{eq:MariasSister}) using the set-up given around equations (9.43) and (9.44) of  
Pace \myand Salvan (1997). To match the notation given there, set
$$\by\equiv \HscrAAi^{-1}\GscrAi\quad 
\mbox{and}\quad\bx\equiv\bU_i^*-\bU_i.$$
Then, in keeping with the displayed equation just before (9.43) of Pace \myand Salvan (1997)
and using their superscript and subscript conventions, we have
$$y^a\equiv\mbox{the $a$th entry of $\by$}\quad\mbox{and}\quad x^a\equiv\mbox{the $a$th entry of $\bx$}.$$
Also,
$$x^{rs}\equiv x^rx^s=\mbox{the $(r,s)$ entry of $\bx\bx^T$}=\mbox{the $(r,s)$ entry 
of $(\bU_i^*-\bU_i)^{\otimes2}$}.$$
Then
$$y^a=x^a+A_{rs}^ax^{rs}+\ldots.$$
where
{\setlength\arraycolsep{1pt}
\begin{eqnarray*}
A_{rs}^a\bx^{rs}&=&\mbox{the $a$th entry of}\ 
\smhalf\HscrAAi^{-1}
\tcbk{\Big[}\HscrdAAAi\bigstar\{(\bU_i^*-\bU_i)(\bU_i^*-\bU_i)^T\}\tcbk{\Big]}\\[2ex]
&=&\mbox{the $a$th entry of}\ 
\smhalf\HscrAAi^{-1}\tcbk{\Big\{}\HscrdAAAi\bigstar(\bx\bx^T)\tcbk{\Big\}}.
\end{eqnarray*}
}
From equations (9.43) and (9.44) of Pace \myand Salvan (1997),
{\setlength\arraycolsep{1pt}
\begin{eqnarray*}
x^a&=&y^a-A_{rs}^a y^{rs}+\ldots\\[1ex]
&=&y^a-\mbox{the $a$th entry of}\ 
\smhalf\HscrAAi^{-1}\tcbk{\Big\{}\HscrdAAAi\bigstar
(\by\by^T)\tcbk{\Big\}}+\ldots\\[1ex]
&=&y^a-\mbox{the $a$th entry of}\ 
\smhalf\HscrAAi^{-1}\tcbk{\Big\{}\HscrdAAAi\bigstar
\tcbk{\Big(}\HscrAAi^{-1}\GscrAi\tcbk{\Big)}^{\otimes2}\tcbk{\Big\}}+\ldots.
\end{eqnarray*}
}
This results in the following three-term approximation of $\bU_i^*$:
\begin{equation}
\bU_i^*=\bU_i+\HscrAAi^{-1}\GscrAi-\smhalf\HscrAAi^{-1}
\Big\{\HscrdAAAi\bigstar\tcbk{\Big(}\HscrAAi^{-1}\GscrAi\GscrAi^T\HscrAAi^{-1}\tcbk{\Big)}\Big\}
+\myOP(n^{-3/2})\bone_{\dR}.
\label{eq:DoncasterPub}
\end{equation}

\section{Score Asymptotic Approximation}\label{sec:scoreApprox}

We are now ready to obtain approximations of the scores $\bSAi$, $\bSBi$ 
and $\bSCi$ with accuracies that are sufficient for the two-term
asymptotic covariance matrices of (\ref{eq:mainResult}).

\subsection{Approximation of $\bSAi$}

For each $1\le r\le\dR$, let $\be_r$ denote the $\dR\times1$ 
vector having $r$th entry equal to 1 and zeroes elsewhere. 

\subsubsection{The (\ref{eq:Miyata}) First Term Contribution}\label{eq:Zebedee}

For each $1\le r\le\dR$, the contribution to the $r$th entry of $\bSAi$ 
from the first term on the right-hand side of (\ref{eq:Miyata}) is
$\mbox{the $r$th entry of}\ \bSigma^{-1}\bU_i^*.$
In view of (\ref{eq:DoncasterPub}) we obtain the following contribution
to $\bSAi$:
$$\bSigma^{-1}\bU_i+\bSigma^{-1}\HscrAAi^{-1}\GscrAi-\smhalf\bSigma^{-1}\HscrAAi^{-1}
\Big\{\HscrdAAAi\bigstar\Big(\HscrAAi^{-1}\GscrAi\GscrAi^T\HscrAAi^{-1}\Big)\Big\}
+\myOP(n^{-3/2})\bone_{\dR}.
$$

\subsubsection{The (\ref{eq:Miyata}) Second Term Contribution}\label{eq:BloodSimple}

Noting that 
$$\nabla\{\be_r^T\bgiA(\bu)\}=\be_r^T\bSigma^{-1}\quad\mbox{and}\quad
\nabla\cS(\bu)=-\cS(\bu)\bSigma^{-1}\bu\tcbk{,}$$
the contribution to the $r$th entry of $\bSAi$ from the second term on the right-hand side 
of (\ref{eq:Miyata}) is 
\begin{equation}
-\phi\be_r^T\bSigma^{-1}
\left\{\dispsumjni b''\big(\linXABUistar\big)\bXAij\bXAij^T\right\}^{-1} 
\bSigma^{-1}\bU_i^*.
\label{eq:CartwheelsOfBulli}
\end{equation}
Substitution of (\ref{eq:DoncasterPub}) into (\ref{eq:CartwheelsOfBulli}) then
leads to the following contribution to $\bSAi$:
$$-\phi\bSigma^{-1}\HscrAAi^{-1}\bSigma^{-1}\bU_i+\myOP(n^{-3/2})\bone_{\dR}.$$

\subsubsection{The (\ref{eq:Miyata}) Third Term Contribution}

Noting that $\nabla^2\{\be_r^T\bgiA(\bu)\}=\bO$, the contribution to $\bSAi$ from 
the third term on the right-hand side of (\ref{eq:Miyata}) is $\bzero$.

\subsubsection{The (\ref{eq:Miyata}) Fourth Term Contribution}\label{sec:CaptainPugwash}

Via arguments similar to those given in Section \ref{eq:BloodSimple}, 
the contribution to $\bSAi$ from the fourth term on the right-hand side of (\ref{eq:Miyata}) is 
$$
-\frac{\tcbk{\phi}}{2}\bSigma^{-1}\HscrAAi^{-1}\Big(\HscrdAAAi\bigstar\HscrAAi^{-1}\Big)
+\myOP(n^{-3/2})\bone_{\dR}.
$$

\subsubsection{The Resultant Score Approximation}

On combining the results of Sections \ref{eq:Zebedee}--\ref{sec:CaptainPugwash}, we obtain
\begin{equation}
{\setlength\arraycolsep{1pt}
\begin{array}{rcl}
  \bSAi&=&\bSigma^{-1}\bU_i+\bSigma^{-1}\HscrAAi^{-1}\GscrAi-{\displaystyle\frac{1}{2}}
           \bSigma^{-1}\HscrAAi^{-1}
\Big\{\HscrdAAAi\bigstar\Big(\HscrAAi^{-1}\GscrAi\GscrAi^T\HscrAAi^{-1}\Big)\Big\}\\[1.5ex]
&&\quad-\phi\bSigma^{-1}\HscrAAi^{-1}\bSigma^{-1}\bU_i
-{\displaystyle\frac{\tcbk{\phi}}{2}}\bSigma^{-1}\HscrAAi^{-1}\Big(\HscrdAAAi\bigstar\HscrAAi^{-1}\Big)
+\myOP(n^{-3/2})\bone_{\dR}.
\end{array}
}
\label{eq:ToomaNotCooma}
\end{equation}

\subsection{Approximation of $\bSBi$}

For each $1\le r\le\dR$, let $\be_r$ denote the $\dB\times1$ 
vector having $r$th entry equal to 1 and zeroes elsewhere. 

\subsubsection{The (\ref{eq:Miyata}) First Term Contribution}

The contribution to $\bSBi$ from the first term on the right-hand side of 
(\ref{eq:Miyata}) is
\begin{equation}
\bgiB(\bU_i^*)=\frac{1}{\phi}\dispsumjni \bXBij\big\{Y_{ij}-b'\big(\linXABUistar\big)\big\}.
\label{eq:BigNoiseFromWinnetka}
\end{equation}
Next note that, with (\ref{eq:DoncasterPub}) as a basis,
{\setlength\arraycolsep{1pt}
\begin{eqnarray*}
&&b'\big(\linXABUistar\big)\\[1ex]
&&\qquad =b'\big((\bbetaA+\bU_i)^T\bXAij+\bbetaB^T\bXBij\big)\\[1ex]
&&\qquad\quad
+\bXAij^T(\bU_i^*-\bU_i)b''\big((\bbetaA+\bU_i)^T\bXAij+\bbetaB^T\bXBij\big)\\[1ex]
&&\qquad\quad+\smhalfdisp\bXAij^T(\bU_i^*-\bU_i)^{\otimes2}\bXAij
b'''\big((\bbetaA+\bU_i)^T\bXAij+\bbetaB^T\bXBij\big)
+\myOP(n^{-3/2})\\[1ex]
&&\qquad =b'\big((\bbetaA+\bU_i)^T\bXAij+\bbetaB^T\bXBij\big)\\[1ex]
&&\qquad\quad+\bXAij^T\left[\HscrAAi^{-1}\GscrAi-\smhalfdisp\HscrAAi^{-1}
\Big\{\HscrdAAAi\bigstar\tcbk{\Big(}\HscrAAi^{-1}\GscrAi\GscrAi^T\HscrAAi^{-1}\tcbk{\Big)}\Big\}\right]\\
&&\qquad\qquad\qquad\times b''\big((\bbetaA+\bU_i)^T\bXAij+\bbetaB^T\bXBij\big)\\[1ex]
&&\qquad\quad+\smhalfdisp\bXAij^T\tcbk{\Big(}\HscrAAi^{-1}\GscrAi\tcbk{\Big)}^{\otimes2}\bXAij
b'''\big((\bbetaA+\bU_i)^T\bXAij+\bbetaB^T\bXBij\big)
+\myOP(n^{-3/2}).
\end{eqnarray*}
}
Substitution of this result into (\ref{eq:BigNoiseFromWinnetka}) leads to the first term of $\bSBi$ 
equalling
{\setlength\arraycolsep{1pt}
\begin{eqnarray*}
&&\frac{1}{\phi}\Big(\GscrBi-\HscrABi^T\HscrAAi^{-1}\GscrAi\Big)
+\frac{1}{2\phi}\HscrABi^T\HscrAAi^{-1}
\Big\{\HscrdAAAi\bigstar\Big(\HscrAAi^{-1}\GscrAi\GscrAi^T\HscrAAi^{-1}\tcbk{\Big)}\Big\}\\[1ex]
&&\qquad\qquad-\frac{1}{2\phi}\Big\{\HscrdAABi\bigstar\Big(\HscrAAi^{-1}\GscrAi\GscrAi^T\HscrAAi^{-1}\tcbk{\Big)}\Big\}
+\myOP(n^{-1/2})\bone_{\dR}.
\end{eqnarray*}
}

\subsubsection{The (\ref{eq:Miyata}) Second Term Contribution}

Noting that 
\begin{equation}
\nabla\{\be_r^T\bgiB(\bu)\}=
-\frac{1}{\phi}\dispsumjni  b''\big(\linXABu\big)\be_r^T\bXAij\bXBij^T
\label{eq:FinleyRockBands}
\end{equation}
and recalling that $\nabla\cS(\bu)=-\cS(\bu)\bSigma^{-1}\bu$, the contribution 
from the second term on the right-hand side of (\ref{eq:Miyata}) to $\bSBi$ is
{\setlength\arraycolsep{1pt}
\begin{eqnarray*}
&&\left\{\dispsumjni b''\big(\linXABUistar\big)\bXAij\bXBij^T\right\}^T\\[0ex]
&&\qquad\times\left\{\dispsumjni b''\big(\linXABUistar\big)\bXAij\bXAij^T\right\}^{-1} 
\bSigma^{-1}\bU_i^*.
\end{eqnarray*}
}
Substitution of (\ref{eq:DoncasterPub}) then leads to the contribution to $\bSBi$
from the second term of (\ref{eq:Miyata}) equalling
$$\HscrABi\tcbk{^T}\HscrAAi^{-1}\bSigma^{-1}\bU_i+\myOP(n^{-1/2})\bone_{\dR}.$$

\subsubsection{The (\ref{eq:Miyata}) Third Term Contribution}

The $r$th entry of the contribution to $\bSBi$ from the third term of (\ref{eq:Miyata}) is
{\setlength\arraycolsep{1pt}
\begin{eqnarray*}
&&\frac{1}{2n}\sum_{s=1}^{\dB}\sum_{t=1}^{\dB}\big[\nabla^2\{\be_r^T\bgiB(\bU_i^*)\}\big]_{st}
\big[\{\nabla^2 h_i(\bU_i^*)\}^{-1}\big]_{st}\\[1ex]
&&\qquad\qquad=\frac{\phi}{2}\sum_{s=1}^{\dB}\sum_{t=1}^{\dB}
\big[\nabla^2\{\be_r^T\bgiB(\bU_i)\}\big]_{st}\big(\HscrAAi^{-1}\big)_{st}+\myOP(n^{-1/2}).
\end{eqnarray*}
}
However, the $(s,t)$ entry of $\nabla^2\{\be_r^T\bgiB(\bU_i)\}$ is
$$-\frac{1}{\phi}\dispsumjni b'''\big(\linXABUi\big)(\be_r^T\bXAij)(\be_s^T\bXAij)(\be_t^T\bXBij)
=-\frac{1}{\phi}\big(\HscrdAABi\big)_{rst}.
$$
Noting (\ref{eq:BlackStarDefn}), the contribution to $\bSBi$ from the third term of (\ref{eq:Miyata}) is
$$-\smhalfdisp\HscrdAABi\bigstar\HscrAAi^{-1}+\myOP(n^{-1/2})\bone_{\dR}.$$

\subsubsection{The (\ref{eq:Miyata}) Fourth Term Contribution}

With the aid of (\ref{eq:FinleyRockBands}), the contribution to $\bSBi$ from the fourth term 
of (\ref{eq:Miyata}) is
{\setlength\arraycolsep{1pt}
\begin{eqnarray*}
&&\frac{1}{2\tcbk{n\phi}}\left\{\dispsumjni  b''\big(\linXABUistar\big)\bXAij\bXBij^T\right\}^{\tcbk{T}}
\{\nabla^2 h_i(\bU_i^*)\}^{-1}
\left[\nabla^3 h_i(\bU_i^*)\bigstar\{\nabla^2 h_i(\bU_i^*)\}^{-1}\right]\\[1ex]
&&\qquad\qquad=\frac{1}{2}\HscrABi^{\tcbk{T}}\HscrAAi^{-1}\Big(\HscrdAAAi\bigstar\HscrAAi^{-1}\Big)
+\myOP(n^{-1/2})\bone_{\dR}.
\end{eqnarray*}
}

\subsubsection{The Resultant Score Approximation}

On combining each of the contributions, we obtain
\begin{equation}
{\setlength\arraycolsep{1pt}
\begin{array}{rcl}
\bSBi&=&{\displaystyle\frac{1}{\phi}}\Big(\GscrBi-\HscrABi^T\HscrAAi^{-1}\GscrAi\Big)
+{\displaystyle\frac{1}{2\phi}}\HscrABi^T\HscrAAi^{-1}
\Big\{\HscrdAAAi\bigstar\Big(\HscrAAi^{-1}\GscrAi\GscrAi^T\HscrAAi^{-1}\tcbk{\Big)}\Big\}\\[1ex]
&&\qquad-{\displaystyle\frac{1}{2\phi}}
\Big\{\HscrdAABi\bigstar\Big(\HscrAAi^{-1}\GscrAi\GscrAi^T\HscrAAi^{-1}\tcbk{\Big)}\Big\}
+\HscrABi^T\HscrAAi^{-1}\bSigma^{-1}\bU_i\\[1.5ex]
&&\qquad-{\displaystyle\frac{1}{2}}\HscrdAABi\bigstar\HscrAAi^{-1}
+{\displaystyle\frac{1}{2}}\HscrABi^T\HscrAAi^{-1}\Big(\HscrdAAAi\bigstar\HscrAAi^{-1}\Big)
+\myOP(n^{-1/2})\bone_{\dR}.
\end{array}
}
\label{eq:CygnetSquad}
\end{equation}

\subsection{Approximation of $\bSCi$}

For each $1\le r\le\smhalf\dR(\dR+1)$ let $\be_r$ denote the $\dR(\dR+1)/2\times1$ vector 
with $1$ in the $r$th position and zeroes elsewhere.

\subsubsection{The (\ref{eq:Miyata}) First Term Contribution}

For each $1\le r\le\smhalf\dR(\dR+1)$, the $r$th entry of the contribution to 
$\bSCi$ from the first term of (\ref{eq:Miyata}) is 
$$\be_r^T\bgiC(\bU_i^*)=
{\displaystyle\frac{1}{2}}\be_r^T\bD_{\dR}^T(\bSigma^{-1}\otimes\bSigma^{-1})
\vecof\big((\bU_i^*)^{\otimes2}\big).$$
Since
{\setlength\arraycolsep{1pt}
\begin{eqnarray*}
(\bU_i^*)^{\otimes2}&=&\tcbk{\Big[}\bU_i+\HscrAAi^{-1}\GscrAi\tcbk{-\smhalf\HscrAAi^{-1}
\Big\{\HscrdAAAi\bigstar\Big(\HscrAAi^{-1}\GscrAi\GscrAi^T\HscrAAi^{-1}\Big)\Big\}}+\myOP(n^{-\tcbk{3/2}})\bone_{\dR}\tcbk{\Big]}^{\otimes2}\\[1ex]
&=&\bU_i^{\otimes2}+\bU_i\GscrAi^T\HscrAAi^{-1}+\HscrAAi^{-1}\GscrAi\bU_i^T\tcbk{+\HscrAAi^{-1}\GscrAi\GscrAi^T\HscrAAi^{-1}}\\[1ex]
&&\tcbk{-\smhalf\bU_i\Big\{\HscrdAAAi\bigstar\Big(\HscrAAi^{-1}
\GscrAi\GscrAi^T\HscrAAi^{-1}\Big)\Big\}^T\HscrAAi^{-1}}\\[1ex]
&&\tcbk{-\smhalf\HscrAAi^{-1}\Big\{\HscrdAAAi\bigstar\Big(\HscrAAi^{-1}
\GscrAi\GscrAi^T\HscrAAi^{-1}\Big)\Big\}\bU_i^T}
+\myOP(n^{-\tcbk{3/2}})\bone_{\dR}^{\otimes2},
\end{eqnarray*}
}
and noting (\ref{eq:MandNthree}) and (\ref{eq:vecABC}),
the contribution to $\bSCi$ from the first term of (\ref{eq:Miyata}) is 
{\setlength\arraycolsep{1pt}
\begin{eqnarray*}
&&\frac{1}{2}\bD_{\dR}^T
\vecof\Big(\bSigma^{-1}\tcbk{\Big[}\bU_i\bU_i^T+2\HscrAAi^{-1}\GscrAi\bU_i^T\tcbk{+\HscrAAi^{-1}\GscrAi\GscrAi^T\HscrAAi^{-1}}\\[0ex]
&&\qquad\qquad\qquad\qquad\tcbk{-\HscrAAi^{-1}\Big\{\HscrdAAAi\bigstar\Big(\HscrAAi^{-1}
\GscrAi\GscrAi^T\HscrAAi^{-1}\Big)\Big\}\bU_i^T}\tcbk{\Big]}\bSigma^{-1}\Big)
+\myOP(n^{-\tcbk{3/2}})\bone_{\dR(\dR+1)/2}.
\end{eqnarray*}
}

\subsubsection{The (\ref{eq:Miyata}) Second Term Contribution}

Noting that, for each $1\le r\le\dR(\dR+1)/2$,
\begin{equation}
\big[\nabla\{\be_r^T\bgiC(\bu)\}\big]^T
=\smhalfdisp\be_r^T\bD_{\dR}^T(\bSigma^{-1}\otimes\bSigma^{-1})
\big\{(\bu\otimes\bI)+(\bI\otimes\bu)\big\}
\label{eq:VenusFlyTrap}
\end{equation}
and keeping in mind that $\nabla\cS(\bu)=-\cS(\bu)\bSigma^{-1}\bu$, the contribution 
from the second term on the right-hand side of (\ref{eq:Miyata}) to $\bSCi$ is
{\setlength\arraycolsep{1pt}
\begin{eqnarray*}
&&-\frac{\phi}{2}
\bD_{\dR}^T(\bSigma^{-1}\otimes\bSigma^{-1})\big\{(\bU_i^*\otimes\bI)+(\bI\otimes\bU_i^*)\big\}\\[0ex]
&&\qquad\qquad\times
\left\{\dispsumjni b''\big(\linXABUistar\big)\bXAij\bXAij^T\right\}^{-1}\bSigma^{-1}\bU_i^*\\[1ex]
&&\qquad=-\frac{\phi}{2}
\bD_{\dR}^T(\bSigma^{-1}\otimes\bSigma^{-1})\big\{(\bU_i\otimes\bI)+(\bI\otimes\bU_i)\big\}
\HscrAAi^{-1}\bSigma^{-1}\bU_i+\myOP(n^{-3/2})\bone_{\dR(\dR+1)/2}\\[1ex]
&&\qquad=-\phi\bD_{\dR}^T\vecof\Big(\bSigma^{-1}\HscrAAi^{-1}\bSigma^{-1}\bU_i\bU_i^T\bSigma^{-1}\Big)
+\myOP(n^{-3/2})\bone_{\dR(\dR+1)/2}.
\end{eqnarray*}
}
The last step makes use of (\ref{eq:MandNone}), (\ref{eq:MandNthree}) and (\ref{eq:vecABC}).

\subsubsection{The (\ref{eq:Miyata}) Third Term Contribution}

The derivation of the (\ref{eq:Miyata}) third term contribution to $\bSCi$ 
benefits from notation and a result concerning the inverse of the $\vecof$
operator. For $d\in\naturalNumbers$, if $\bb$ is a $d^2\times1$ vector then
$\vecof^{-1}(\bb)$ is the $d\times d$ matrix such that 
$\vecof\big(\vecof^{-1}(\bb)\big)=\bb$. 

\noindent
\begin{lemma}
Let $d\in\naturalNumbers$, $\ba$ be a $d\times1$ vector and 
$\bb$ be a $d^2\times 1$ vector. Then 
$$(\ba^T\otimes\bI)\bb=\vecof^{-1}(\bb)\ba\quad\mbox{and}\quad
(\bI\otimes\ba^T)\bb=\vecof^{-1}(\bb)^T\ba.
$$
\label{lem:vecInvRes}
\end{lemma}

\noindent
Lemma \ref{lem:vecInvRes} is a relatively simple consequence of 
(\ref{eq:vecABC}). To prove the first part of Lemma \ref{lem:vecInvRes}, 
note that its right-hand side is
$$\vecof^{-1}(\bb)\ba=\vecof\big(\vecof^{-1}(\bb)\ba\big)
=\vecof\big(\bI\vecof^{-1}(\bb)\ba\big)=
(\ba^T\otimes\bI)\vecof\big(\vecof^{-1}(\bb)\big)
=(\ba^T\otimes\bI)\bb.
$$
The proof of the second part of Lemma \ref{lem:vecInvRes} is similar.
\vskip10mm

For each $1\le r\le\smhalf\dR(\dR+1)$, the $r$th entry of the contribution to 
$\bSCi$ from the third term of (\ref{eq:Miyata}) is 
$$\frac{1}{2n}\tr\left[\{\nabla^2 h_i(\bU_i)\}^{-1}\nabla^2\{\be_r^T\bgiC(\bU_i^*)\}\right].$$
Next note \tcbk{from \eqref{eq:VenusFlyTrap}} that 
\begin{equation*}
d\{\be_r^T\bgiC(\bu)\}=\smhalfdisp\be_r^T\bD_{\dR}^T(\bSigma^{-1}\otimes\bSigma^{-1})
\big\{(\bu\otimes\bI)+(\bI\otimes\bu)\big\}d\bu.
\label{eq:DonohoIceCream}
\end{equation*}
Using Lemma \ref{lem:vecInvRes} we then have
{\setlength\arraycolsep{1pt}
\begin{eqnarray*}
2d^2\{\be_r^T\bgiC(\bu)\}&=&\be_r^T\bD_{\dR}^T(\bSigma^{-1}\otimes\bSigma^{-1})
\big\{(d\bu\otimes\bI)+(\bI\otimes d\bu)\big\}d\bu\\[1ex]
&=&\Big[\big\{(\bSigma^{-1}\otimes\bSigma^{-1})\bD_{\dR}\be_r\big\}^T(d\bu\otimes\bI)
+\big\{(\bSigma^{-1}\otimes\bSigma^{-1})\bD_{\dR}\be_r\big\}^T(\bI\otimes d\bu)\Big]d\bu\\[1ex]
&=&(d\bu)^T\Big[\vecof^{-1}\Big((\bSigma^{-1}\otimes\bSigma^{-1})\bD_{\dR}\be_r\Big)
+\vecof^{-1}\Big((\bSigma^{-1}\otimes\bSigma^{-1})\bD_{\dR}\be_r\Big)^T
\Big]d\bu.
\end{eqnarray*}
}
From the second identification theorem of matrix differential calculus
(e.g.\ Magnus \myand Neudecker, 1999) we then have
$$\nabla^2\{\be_r^T\bgiC(\bu)\}=\smhalfdisp\vecof^{-1}\Big((\bSigma^{-1}\otimes\bSigma^{-1})\bD_{\dR}\be_r\Big)
+\smhalfdisp\vecof^{-1}\Big((\bSigma^{-1}\otimes\bSigma^{-1})\bD_{\dR}\be_r\Big)^T
$$
which does not depend on $\bu$. Therefore $\nabla^2\bg_k(\bU_i^*)$ is a symmetric matrix
that depends only on $\bSigma$, which we denote as follows:
$$\bQ(\bSigma;r)\equiv \smhalfdisp\vecof^{-1}\Big((\bSigma^{-1}\otimes\bSigma^{-1})\bD_{\dR}\be_r\Big)
+\smhalfdisp\vecof^{-1}\Big((\bSigma^{-1}\otimes\bSigma^{-1})\bD_{\dR}\be_r\Big)^T.
$$

Now note that 
{\setlength\arraycolsep{1pt}
\begin{eqnarray*}
\frac{1}{2n}\tr\left[\{\nabla^2 h_i(\bU_i^*)\}^{-1}\nabla^2\{\be_r^T\bgiC(\bU_i^*)\}\right]
&=&\tr\left[\{\nabla^2 h(_i\bU_i^*)\}^{-1}\bQ(\bSigma;r)\right]\\[1ex]
&=&\frac{1}{2n}\tr\left[\{\nabla^2 h_i(\bU_i)\}^{-1}\bQ(\bSigma;r)\right]+\myOP(n^{-3/2})\\[1ex]
&=&\frac{\phi}{2}\tr\left\{\HscrAAi^{-1}\bQ(\bSigma;r)\right\}+\myOP(n^{-3/2}).
\end{eqnarray*}
}
The $r$th entry of the leading term of the contribution to $\bSCi$ from
the third term on the right-hand side of (\ref{eq:Miyata}) is
{\setlength\arraycolsep{1pt}
\begin{eqnarray*}
&&\frac{\phi}{2}\tr\left\{\HscrAAi^{-1}\bQ(\bSigma;r)\right\}\\
&&\qquad\qquad=\frac{\phi}{4}\vecof\big(\HscrAAi^{-1}\big)^T
\vecof\Big(\vecof^{-1}\Big((\bSigma^{-1}\otimes\bSigma^{-1})\bD_{\dR}\be_r\Big)\Big)\\[1ex]
&&\qquad\qquad\qquad+\frac{\phi}{4}\vecof\big(\HscrAAi^{-1}\big)^T
\vecof\Big(\vecof^{-1}\Big((\bSigma^{-1}\otimes\bSigma^{-1})\bD_{\dR}\be_r\Big)^T\Big)\\[1ex]
&&\qquad\qquad=\frac{\phi}{4}\vecof\big(\HscrAAi^{-1}\big)^T
\vecof\Big(\vecof^{-1}\Big((\bSigma^{-1}\otimes\bSigma^{-1})\bD_{\dR}\be_r\Big)\Big)\\[1ex]
&&\qquad\qquad\qquad+\frac{\phi}{4}\vecof\big(\HscrAAi^{-1}\big)^T
\bK_{\dR}\vecof\Big(\vecof^{-1}\Big((\bSigma^{-1}\otimes\bSigma^{-1})\bD_{\dR}\be_r\Big)\Big)\\[1ex]
&&\qquad\qquad
=\frac{\phi}{4}\vecof\big(\HscrAAi^{-1}\big)^T(\bSigma^{-1}\otimes\bSigma^{-1})\bD_{\dR}\be_r
+\frac{\phi}{4}\vecof\big(\HscrAAi^{-1}\big)^T\bK_{\dR}(\bSigma^{-1}\otimes\bSigma^{-1})\bD_{\dR}\be_r\\[1ex]
&&\qquad\qquad
=\frac{\phi}{4}\vecof\big(\HscrAAi^{-1}\big)^T(\bSigma^{-1}\otimes\bSigma^{-1})\bD_{\dR}\be_r
+\frac{\phi}{4}\vecof\big(\HscrAAi^{-1}\big)^T(\bSigma^{-1}\otimes\bSigma^{-1})\bK_{\dR}\bD_{\dR}\be_r\\[1ex]
&&\qquad\qquad
=\frac{\phi}{4}\vecof\big(\HscrAAi^{-1}\big)^T(\bSigma^{-1}\otimes\bSigma^{-1})\bD_{\dR}\be_r
+\frac{\phi}{4}\vecof\big(\HscrAAi^{-1}\big)^T(\bSigma^{-1}\otimes\bSigma^{-1})\bD_{\dR}\be_r\\[1ex]
&&\qquad\qquad
=\frac{\phi}{2}\be_r^T\bD_{\dR}^T(\bSigma^{-1}\otimes\bSigma^{-1})\vecof\big(\HscrAAi^{-1}\big)
=\frac{\phi}{2}\be_r^T\bD_{\dR}^T\vecof\big(\bSigma^{-1}\HscrAAi^{-1}\bSigma^{-1}\big).
\end{eqnarray*}
}
Hence, contribution to $\bSCi$ from the third term on the right-hand side of (\ref{eq:Miyata}) is
$$\frac{\phi}{2}\bD_{\dR}^T\vecof\big(\bSigma^{-1}\HscrAAi^{-1}\bSigma^{-1}\big)
+\myOP(n^{-3/2})\bone_{\dR(\dR+1)/2}.$$

\subsubsection{The (\ref{eq:Miyata}) Fourth Term Contribution}

For each $1\le r\le\smhalf\dR(\dR+1)$, the $r$th entry of the contribution 
to $\bSCi$ from the fourth term on the right-hand side of (\ref{eq:Miyata}) is
$$-\frac{1}{2n}\big[\nabla\{\be_r^T\bgiC(\bU_i^*)\}\big]^T\{\nabla^2 h_i(\bU_i^*)\}^{-1}\big[\nabla^3h_i(\bU_i^*)
\bigstar \{\nabla^2 h_i(\bU_i^*)\}^{-1}\big].$$
Noting (\ref{eq:VenusFlyTrap}) \tcbk{and using \eqref{eq:ACkronBD}}, it follows that
the contribution to $\bSCi$ from the fourth term on the right-hand side 
of (\ref{eq:Miyata}) is
{\setlength\arraycolsep{1pt}
\begin{eqnarray*}
&&-\frac{1}{4n}\bD_{\dR}^T(\bSigma^{-1}\otimes\bSigma^{-1})
\big\{(\bU_i^*\otimes\bI)+(\bI\otimes\bU_i^*)\big\}
\{\nabla^2 h_i(\bU_i^*)\}^{-1}\big[\nabla^3h_i(\bU_i^*)
\bigstar \{\nabla^2 h_i(\bU_i^*)\}^{-1}\big]\\[2ex]
&&\ =-\frac{1}{4n}\bD_{\dR}^T
\big[\{(\bSigma^{-1}\bU_i^*)\otimes\bSigma^{-1}\}+\{\bSigma^{-1}\otimes(\bSigma^{-1}\bU_i^*)\}\big]\\
&&\qquad\qquad\qquad\qquad\times\{\nabla^2 h_i(\bU_i^*)\}^{-1}\big[\nabla^3h_i(\bU_i^*)
\bigstar \{\nabla^2 h_i(\bU_i^*)\}^{-1}\big]\\[2ex]
&&\ =-\frac{\phi}{4}\bD_{\dR}^T
\big[\{(\bSigma^{-1}\bU_i)\otimes\bSigma^{-1}\}+\{\bSigma^{-1}\otimes(\bSigma^{-1}\bU_i)\}\big]
\HscrAAi^{-1}\tcbk{\Big(}\HscrdAAAi\bigstar\HscrAAi^{-1}\tcbk{\Big)}+\myOP(n^{-3/2})\bone_{\dR(\dR+1)/2}
\\[2ex]
&&\ \tcbk{=-\frac{\phi}{2}
\bD_{\dR}^T\vecof\Big(\bSigma^{-1}\Big\{\HscrAAi^{-1}\tcbk{\Big(}\HscrdAAAi\bigstar\HscrAAi^{-1}\tcbk{\Big)}\bU_i^T
\Big\}\bSigma^{-1}\Big)+\myOP(n^{-3/2})\bone_{\dR(\dR+1)/2}}
\end{eqnarray*}
}
\tcbk{where the last step follows from application of (\ref{eq:MandNthree}) and (\ref{eq:vecABC}).}

\subsubsection{The Resultant Score Approximation}

The resultant approximation of $\bSCi$ is 
\begin{equation}
{\setlength\arraycolsep{1pt}
\begin{array}{rcl}
\bSCi&=&\smhalfdisp\bD_{\dR}^T\vecof\Big(\bSigma^{-1}\Big[\bU_i\bU_i^T-\bSigma
+2\HscrAAi^{-1}\GscrAi\bU_i^T+\HscrAAi^{-1}\GscrAi\GscrAi^T\HscrAAi^{-1}\\[1ex]
&&\qquad\qquad +\phi\HscrAAi^{-1}-2\phi\HscrAAi^{-1}\bSigma^{-1}\bU_i\bU_i^T\\[1ex]
&&\qquad\qquad-\HscrAAi^{-1}\Big\{\HscrdAAAi\bigstar\tcbk{\Big(}\HscrAAi^{-1}
\GscrAi\GscrAi^T\HscrAAi^{-1}\tcbk{\Big)}\Big\}\bU_i^T\\[1ex]
&&\qquad\qquad -\phi\HscrAAi^{-1}\tcbk{\Big(}\HscrdAAAi\bigstar\HscrAAi^{-1}\tcbk{\Big)}\bU_i^T
\Big]\bSigma^{-1}\Big)+\myOP(n^{-3/2})\bone_{\dR(\dR+1)/2}.
\end{array}
}
\label{eq:MexicanHat}
\end{equation}

\section{Score Outer Product Conditional Moments Approximation}\label{sec:outProd}

The $i$th term of the Fisher information matrix of $\big(\bbeta,\vech(\bSigma)\big)$
is a $3\times3$ block partitioned matrix with the blocks corresponding
to the various moments of pairwise outer products, conditional
on $\bX_i$. The relevant approximations involve repeated use of (\ref{eq:Ulverstone}) 
and and keeping track of orders of magnitude.

\subsection{Approximation of $E(\bSAi^{\otimes2}|\bX_i)$}

Using (\ref{eq:ToomaNotCooma}), (\ref{eq:Ulverstone}) and 
standard algebraic steps we have
\begin{equation}
{\setlength\arraycolsep{1pt}
\begin{array}{rcl}
&&E(\bS_{Ai}^{\otimes 2}|\bX_i)
=\bSigma^{-1}\\[1ex]
&&\qquad+\phi\bSigma^{-1}E\Big(\HscrAAi^{-1}-\bU_i\bU_i^T\bSigma^{-1}\HscrAAi^{-1}         
-\HscrAAi^{-1}\bSigma^{-1}\bU_i\bU_i^T\
|\bX_i\Big)\bSigma^{-1}\\[1ex]
&&\qquad-\phi\bSigma^{-1}E\Big\{\HscrAAi^{-1}
\tcbk{\Big(}\HscrdAAAi\bigstar\HscrAAi^{-1}\tcbk{\Big)}\bU_i^T+
\bU_i\tcbk{\Big(}\HscrdAAAi\bigstar\HscrAAi^{-1}\tcbk{\Big)}^T\HscrAAi^{-1}
\Big|\bX_i\Big\}\bSigma^{-1}\\[1ex]
&&\qquad+\myOP(n^{-2})\bone^{\otimes2}_{\dR}.
\end{array}
}
\label{eq:AmakiSauce}
\end{equation}

\subsection{Approximation of $E(\bSBi^{\otimes2}|\bX_i)$}

From (\ref{eq:CygnetSquad}) and (\ref{eq:Ulverstone}) we obtain
\begin{equation}
E(\bS_{Bi}^{\otimes2}|\bX_i)= \frac{1}{\phi}E\Big(\HscrBBi-\HscrABi^T\HscrAAi^{-1}\HscrABi\big|\bX_i\Big)
+\myOP(1)\bone_{\dB}^{\otimes2}.
\label{eq:HappyRugbyLeague}
\end{equation}

\subsection{Approximation of $E(\bSCi^{\otimes2}|\bX_i)$}

After some long-winded, but relatively straightforward, matrix algebra that
involves application of (\ref{eq:Ulverstone}) we have \tcbk{from \eqref{eq:MexicanHat} that}
\begin{equation}
{\setlength\arraycolsep{1pt}
\begin{array}{rcl}
&&E(\bSCi^{\otimes2}|\bX_i)=\smhalfdisp\bD_{\dR}^T(\bSigma^{-1}\otimes\bSigma^{-1})\bD_{\dR}
+\frac{\phi}{2}\bD_{\dR}^T(\bSigma^{-1}\otimes\bSigma^{-1})E\Big[2(\bU_i\bU_i^T)\otimes\HscrAAi^{-1} \\[1ex]
&&\quad+\vecof(\bU_i\bU_i^T-\bSigma)\vecof\left(
\HscrAAi^{-1}\bSigma^{-1}\Big\{\bSigma- \bU_i\bU_i^T-\bSigma 
\tcbk{\Big(}\HscrdAAAi\bigstar\HscrAAi^{-1}\tcbk{\Big)}\bU_i^T\Big\}
\right)^T\\[1ex]
&&\quad+
\vecof\left(
\HscrAAi^{-1}\bSigma^{-1}\Big\{\bSigma- \bU_i\bU_i^T-\bSigma 
\tcbk{\Big(}\HscrdAAAi\bigstar\HscrAAi^{-1}\tcbk{\Big)}\bU_i^T\Big\}
\right)\vecof(\bU_i\bU_i^T-\bSigma)^T\Big|\bX_i\Big]\\[1ex]
&&\quad\times(\bSigma^{-1}\otimes\bSigma^{-1})\bD_{\dR}
+ \myOP(n^{-2})\bone_{\dR(\dR+1)/2}^{\otimes2}.
\end{array}
}
\label{eq:BarkBustersDenver}
\end{equation}

\subsection{Approximation of $E(\bSAi\bSBi^T|\bX_i)$}

Multiplication of (\ref{eq:ToomaNotCooma}) by the transpose of
(\ref{eq:CygnetSquad}), taking expectations conditional on $\bX_i$ and 
use of (\ref{eq:Ulverstone}) leads to 
\begin{equation}
{\setlength\arraycolsep{1pt}
\begin{array}{rcl}
E(\bSAi\bSBi^T|\bX_i)&=&\bSigma^{-1}
                         E\Big\{\bU_i\bU_i^T\bSigma^{-1}\HscrAAi^{-1}\HscrABi
                         -\bU_i\tcbk{\Big(}\HscrdAABi\bigstar\HscrAAi^{-1}\tcbk{\Big)}^T\\[1ex]
&&\qquad +\bU_i\Big(\HscrdAAAi\bigstar\HscrAAi^{-1}\Big)^T\HscrAAi^{-1}\HscrABi\Big|\bX_i\Big\}
+\myOP(n^{-1})\bone_{\dR}\bone_{\dB}^T.
\end{array}
\label{eq:BeenToOmeo}
}
\end{equation}

\subsection{Approximation of $E(\bSAi\bSCi^T|\bX_i)$}

An important aspect of the $E(\bSAi\bSCi^T|\bX_i)$ approximation
is that, even though 
$$\bSAi=\myOP(1)\bone_{\dR}\quad\mbox{and}\quad\bSCi=\myOP(1)\bone_{\dR(\dR+1)/2}$$
we can establish that 
\begin{equation}
E(\bSAi\bSCi^T|\bX_i)=\myOP(n^{-1})\bone_{\dR}\bone_{\dR(\dR+1)/2}^T,
\label{eq:ManFromHongKong}
\end{equation}
which indicates a degree of asymptotic orthogonality between $\bbetaA$
and $\bSigma$. An illustrative cancellation, involving the leading terms 
of each score, is
$$E\{\bSigma^{-1}\bU_i\bSigma^{-1}(\bU_i\bU_i^T-\bSigma)\bSigma^{-1}|\bX_i\}\bD_{\dR}
=\bSigma^{-1}\bU_i\bSigma^{-1}\{E(\bU_i\bU_i^T|\bX_i)-\bSigma\}\bSigma^{-1}\bD_{\dR}
=\bO.
$$
As will be shown in Section \ref{sec:ApproxAsyCov}, approximation (\ref{eq:ManFromHongKong})
is sufficient for (\ref{eq:mainResult}).

\subsection{Approximation of $E(\bSBi\bSCi^T|\bX_i)$}

Multiplication of (\tcbk{\ref{eq:CygnetSquad}}) by the transpose of 
(\ref{eq:MexicanHat}), and similar arguments, leads to
\begin{equation}
{\setlength\arraycolsep{1pt}
\begin{array}{rcl}
&&E\big(\bSBi\bSCi^T\big|\bX_i\big)
=\displaystyle{\frac{1}{2}}E\Big[\Big\{\HscrABi^T\HscrAAi^{-1}\bSigma^{-1}\bU_i
-\HscrdAABi\bigstar\HscrAAi^{-1}\\[1.5ex]
&&\qquad+\HscrABi^T\HscrAAi^{-1}
\Big(\HscrdAAAi\bigstar\HscrAAi^{-1}\Big)\Big\}
\vecof(\bU_i\bU_i^T-\bSigma)^T
\Big|\bX_i\Big](\bSigma^{-1}\otimes\bSigma^{-1})\bD_{\dR}\\[1.5ex]
&&\qquad+\myOP(n^{-1})\bone_{\dB}\bone^T_{\dR(\dR+1)/2}.
\end{array}
}
\label{eq:AndPigsDoFly}
\end{equation}

\section{The Fisher Information Matrix}\label{sec:FishInfo}

The Fisher information matrix of $\big(\bbeta,\vech(\bSigma)\big)$ is
{\setlength\arraycolsep{1pt}
\begin{eqnarray*}
I\big(\bbeta,\vech(\bSigma)\big)
&=&\left[
\begin{array}{cc}
\bM_{11}   & \bM_{12} \\[1ex]
\bM_{12}^T & \bM_{22}
\end{array}
\right]
\quad\mbox{where}\quad
\bM_{11}\equiv \sumim\left[
\begin{array}{cc}
E(\bSAi^{\otimes2}|\bX_i) &  E(\bSAi\bSBi^T|\bX_i) \\[1ex]
E(\bSBi\bSAi^T|\bX_i) &  E(\bSBi^{\otimes2}|\bX_i)  
\end{array}
\right]\tcbk{,}
\end{eqnarray*}
}
$$
\bM_{12}\equiv \sumim\left[
\begin{array}{c}
E(\bSAi\bSCi^T|\bX_i) \\[1ex]
E(\bSBi\bSCi^T|\bX_i)  
\end{array}
\right]
\quad\mbox{and}\quad
\bM_{22}\equiv \sumim E(\bSCi^{\otimes2}|\bX_i).
$$
The results of the previous section lead to high-order asymptotic approximation
of the matrix $I\big(\bbeta,\vech(\bSigma)\big)$. In the next section we show that
inversion of this approximate Fisher information matrix leads to two-term covariance matrix
approximations for the maximum likelihood estimators.

\section{Approximation of Covariance Matrices of Estimators}\label{sec:ApproxAsyCov}

The dominant terms in the approximation of 
$$\Cov\big(\bbetaMLE|\bXsc\big)\quad\mbox{and}\quad\Cov\big(\vech(\bSigmaMLE)|\bXsc\big)$$
correspond to the $\dF\times\dF$ and $\frac{1}{2}\dR(\dR+1)\times\frac{1}{2}\dR(\dR+1)$ 
diagonal blocks of 
$$I\big(\bbeta,\vech(\bSigma)\big)^{-1}.$$
We now treat each of these in turn in the upcoming subsections\tcbk{,
which make extensive use of block matrix inversions}. \tcbk{If a matrix
is partitioned into four blocks $\bA$, $\bB$, $\bC$ and $\bD$, then}
$$\tcbk{\left[
\begin{array}{cc}
\bA & \bB\\[1ex]
\bC & \bD
\end{array}
\right]^{-1}=\left[
\begin{array}{cc}
\bA^{-1}+\bA^{-1}\bB(\bD-\bC\bA^{-1}\bB)^{-1}\bC\bA^{-1} & -\bA^{-1}\bB(\bD-\bC\bA^{-1}\bB)^{-1}\\[1ex]
-(\bD-\bC\bA^{-1}\bB)^{-1}\bC\bA^{-1} & (\bD-\bC\bA^{-1}\bB)^{-1}
\end{array}
\right]}$$
\tcbk{or, equivalently,}
$$\tcbk{\left[
\begin{array}{cc}
\bA & \bB\\[1ex]
\bC & \bD
\end{array}
\right]^{-1}=\left[
\begin{array}{cc}
(\bA-\bB\bD^{-1}\bC)^{-1} & -(\bA-\bB\bD^{-1}\bC)^{-1}\bB\bD^{-1}\\[1ex]
-\bD^{-1}\bC(\bA-\bB\bD^{-1}\bC)^{-1} & \bD^{-1}+\bD^{-1}\bC(\bA-\bB\bD^{-1}\bC)^{-1}\bB\bD^{-1}
\end{array}
\right].}$$
Another result that is repeatedly used in the following subsections is 
$$(\bA-\bB)^{-1}=\sum_{k=0}^{\infty}(\bA^{-1}\bB)^k\bA^{-1}$$ 
for $\bA$ and $\bB$ invertible matrices of the same size and
such that the spectral radius of $\bA^{-1}\bB$ is less than $1$.

\subsection{Two-Term Approximation of $\Cov\big(\bbetaMLE|\bXsc\big)$}

The dominant terms of $\Cov\big(\bbetaMLE|\bXsc\big)$ correspond to
$$\mbox{the upper left $\dF\times\dF$ block of 
$I\big(\bbeta,\vech(\bSigma)\big)^{-1}$}
=\big(\bM_{11}-\bM_{12}\bM_{22}^{-1}\bM_{12}^T\big)^{-1}.
$$
Based on (\ref{eq:AmakiSauce}), (\ref{eq:HappyRugbyLeague}) and (\ref{eq:BeenToOmeo}) we have 
$$\bM_{11}=\left[
\begin{array}{cc}
m\bSigma^{-1}-{\displaystyle\frac{\phi m}{n}}\bSigma^{-1}\bKscAA\bSigma^{-1}
+\myOP(mn^{-2})\bone_{\dR}^{\otimes2} 
& m\bSigma^{-1}\bKscAB+\myOP(mn^{-1})\bone_{\dR}\bone_{\dB}^T\\[1ex]
m\bSigma^{-1}\bKscAB^T+\myOP(mn^{-1})\bone_{\dB}\bone_{\dR}^T 
& {\displaystyle\frac{mn}{\phi}}\bKscBB+\myOP(m)\bone_{\dB}^{\otimes2} 
\end{array}
\right]
$$
where 
{\setlength\arraycolsep{1pt}
\begin{eqnarray*}
\bKscAA&\equiv&\frac{n}{m}\sumim E\Big\{\bU_i\bU_i^T\bSigma^{-1}\HscrAAi^{-1}         
+\HscrAAi^{-1}\bSigma^{-1}\bU_i\bU_i^T-\HscrAAi^{-1}\\[0ex]
&&\qquad\qquad+\HscrAAi^{-1}
\tcbk{\Big(}\HscrdAAAi\bigstar\HscrAAi^{-1}\tcbk{\Big)}\bU_i^T
+\bU_i\tcbk{\Big(}\HscrdAAAi\bigstar\HscrAAi^{-1}\tcbk{\Big)}^T\HscrAAi^{-1}
\Big|\bX_i\Big\},\\[1ex]
\bKscAB&\equiv&\frac{1}{m}\sumim 
                E\Big\{\bU_i\bU_i^T\bSigma^{-1}\HscrAAi^{-1}\HscrABi-\bU_i
                \tcbk{\Big(}\HscrdAABi\bigstar\HscrAAi^{-1}\tcbk{\Big)}^T\\[0ex]
&&\qquad\qquad\qquad\qquad\qquad
+\bU_i\Big(\HscrdAAAi\bigstar\HscrAAi^{-1}\Big)^T\HscrAAi^{-1}\HscrABi\Big|\bX_i\Big\}\\[1ex]
\mbox{and}\quad
\bKscBB&\equiv&\frac{1}{mn}\sumim E\Big(\HscrBBi-\HscrABi^T\HscrAAi^{-1}\HscrABi\Big|\bX_i\Big)
\end{eqnarray*}
}
are matrices with all entries being $\myOP(1)$. As consequences of (\ref{eq:BarkBustersDenver}), 
(\ref{eq:ManFromHongKong}) and (\ref{eq:AndPigsDoFly}) we have
\begin{equation}
\bM_{22}^{-1}=\myOP(m^{-1})\bone_{\dR(\dR+1)/2}^{\otimes2}
\quad\mbox{and}\quad
\bM_{12}=\left[\begin{array}{c}
\myOP(mn^{-1})\bone_{\dR}\bone_{\dR(\dR+1)/2}^T\\[1ex]
\myOP(m)\bone_{\dB}\bone_{\dR(\dR+1)/2}^T
\end{array}
\right].
\label{eq:WidnesFactory}
\end{equation}
Therefore
$$\bM_{12}\bM_{22}^{-1}\bM_{12}^T
=\left[
\begin{array}{cc}
\myOP(mn^{-2})\bone_{\dR}^{\otimes2} & \myOP(mn^{-1})\bone_{\dR}\bone_{\dB}^T\\[1ex]
\myOP(mn^{-1})\bone_{\dB}\bone_{\dR}^T &\quad \myOP(m)\bone_{\dB}^{\otimes2}
\end{array}
\right].
$$
From these results for $\bM_{11}$ and $\bM_{12}\bM_{22}^{-1}\bM_{12}^T$, it 
follows that 
\begin{eqnarray*}
&&\bM_{11}-\bM_{12}\bM_{22}^{-1}\bM_{12}^T=\\[1ex]
&&\qquad\left[
\begin{array}{cc}
m\bSigma^{-1}-{\displaystyle\frac{\phi m}{n}}\bSigma^{-1}\bKscAA\bSigma^{-1}
+\myOP(mn^{-2})\bone_{\dR}^{\otimes2} 
& m\bSigma^{-1}\bKscAB+\myOP(mn^{-1})\bone_{\dR}\bone_{\dB}^T\\[1ex]
m\bSigma^{-1}\bKscAB^T+\myOP(mn^{-1})\bone_{\dB}\bone_{\dR}^T 
& {\displaystyle\frac{mn}{\phi}}\bKscBB+\myOP(m)\bone_{\dB}^{\otimes2} 
\end{array}
\right].
\end{eqnarray*}
The upper left $\dR\times\dR$ block of 
$\big(\bM_{11}-\bM_{12}\bM_{22}^{-1}\bM_{12}^T\big)^{-1}$ is
{\setlength\arraycolsep{1pt}
\begin{eqnarray*}
&&\left\{
m\bSigma^{-1}-{\displaystyle\frac{\phi m}{n}}\bSigma^{-1}\big(\bKscAA+\bKscAB\bKscBB^{-1}\bKscAB^T\big)\bSigma^{-1}
+\myOP(mn^{-2})\bone_{\dR}^{\otimes2}
\right\}^{-1}\\[1ex]
&&
\qquad\qquad=
\frac{1}{m}\Big\{
\bI-{\displaystyle\frac{\phi}{n}}\big(\bKscAA+\bKscAB\bKscBB^{-1}\bKscAB^T\big)\bSigma^{-1}
+\myOP(n^{-2})\bone_{\dR}^{\otimes2}\Big\}^{-1}\bSigma\\[1ex]
&&\qquad\qquad=\frac{\bSigma}{m}+\frac{\phi\big(\bKscAA+\bKscAB\bKscBB^{-1}\bKscAB^T\big)}{mn}
+\myOP(m^{-1}n^{-2})\bone_{\dR}^{\otimes2}.
\end{eqnarray*}
}
The upper right $\dR\times\dB$ block of $\big(\bM_{11}-\bM_{12}\bM_{22}^{-1}\bM_{12}^T\big)^{-1}$ is
{\setlength\arraycolsep{1pt}
\begin{eqnarray*}
&&-\left\{\frac{\bSigma}{m}+\myOP(m^{-1}n^{-1})\bone_{\dR}^{\otimes2}\right\}
\Big\{m\bSigma^{-1}\bKscAB+\myOP(mn^{-1})\bone_{\dR}\bone_{\dB}^T\Big\}
\left\{{\displaystyle\frac{mn}{\phi}}\bKscBB+\myOP(m)\bone_{\dB}^{\otimes2}\right\}^{-1}\\[1ex]
&&\qquad\qquad=-\frac{\phi\bKscAB\bKscBB^{-1}}{mn}+\myOP(m^{-1}n^{-2})\bone_{\dR}\bone_{\dB}^T.
\end{eqnarray*}
}
\noindent
The lower right $\dB\times\dB$ block of $\big(\bM_{11}-\bM_{12}\bM_{22}^{-1}\bM_{12}^T\big)^{-1}$ is
$$\frac{\phi\bKscBB^{-1}}{mn}+\myOP(m^{-1}n^{-2})\bone_{\dB}^{\otimes2}.$$
Therefore, 
{\setlength\arraycolsep{1pt}
\begin{eqnarray*}
\Cov\big(\bbetaMLE|\bXsc\big)&=&\frac{1}{m}\left[
\begin{array}{cc}
\bSigmaZero &\quad \bO \\[1ex]
\bO     &\quad  \bO
\end{array}
\right]
+
\frac{\phi}{mn}
\left[
\begin{array}{cc}
(\bKscAA^0)^{-1} &\qquad (\bKscAA^0)^{-1}\bKscAB^0\\[1ex]
(\bKscAB^0)^T(\bKscAA^0)^{-1}&\qquad \bKscBB^0+(\bKscAB^0)^T
(\bKscAA^0)^{-1}\bKscAB^0
\end{array}
\right]^{-1}\\[1ex]
&&\qquad+\myOP(m^{-1}n^{-2})\bone_{\dF}^{\otimes2}
\end{eqnarray*}
}
where, for example, $\bKscAA^0$ is the $\bKscAA$ quantity with $\bbeta$ set to $\bbetaZero$
and $\bSigma$ set to $\bSigmaZero$.

\subsection{Two-Term Approximation of $\Cov\big(\vech(\bSigmaMLE)|\bXsc\big)$}

The dominant terms of $\Cov\big(\vech(\bSigmaMLE)|\bXsc\big)$ correspond to
\begin{eqnarray*}
&&\mbox{the lower right $\smhalfdisp\dR(\dR+1)\times\smhalfdisp\dR(\dR+1)$ block of 
$I\big(\bbeta,\vech(\bSigma)\big)^{-1}$}\\[1ex]
&&\qquad\qquad\qquad\qquad\qquad\qquad\qquad\qquad\qquad\qquad
=\big(\bM_{22}-\bM_{12}^T\bM_{11}^{-1}\bM_{12}\big)^{-1}.
\end{eqnarray*}
From (\ref{eq:BarkBustersDenver}) 
{\setlength\arraycolsep{1pt}
\begin{eqnarray*}
\sumim E(\bSCi\bSCi^T|\bX_i)
&=&\frac{m}{2}\bD_{\dR}^T(\bSigma^{-1}\otimes\bSigma^{-1})\bD_{\dR}
-\frac{m\phi}{n}\bD_{\dR}^T(\bSigma^{-1}\otimes\bSigma^{-1})\bKscCC(\bSigma^{-1}\otimes\bSigma^{-1})
\bD_{\dR}\\[0ex]
&&\qquad+\myOP(mn^{-2})\bone_{\dR(\dR+1)/2}^{\otimes2}
\end{eqnarray*}
}
with the following $\myOP(1)$ matrix: 
{\setlength\arraycolsep{1pt}
\begin{eqnarray*}
\bKscCC&=&\frac{n}{m}\sumim 
E\Big[\displaystyle{\frac{1}{2}}\vecof(\bSigma-\bU_i\bU_i^T)\vecof\left(
\HscrAAi^{-1}\bSigma^{-1}\Big\{\bSigma- \bU_i\bU_i^T-\bSigma \tcbk{\Big(}\HscrdAAAi\bigstar\HscrAAi^{-1}\tcbk{\Big)}\bU_i^T\Big\}
\right)^T\\[1ex]
&&\qquad\qquad+
\displaystyle{\frac{1}{2}}\vecof\left(
\HscrAAi^{-1}\bSigma^{-1}\Big\{\bSigma- \bU_i\bU_i^T-\bSigma \tcbk{\Big(}\HscrdAAAi\bigstar\HscrAAi^{-1}\tcbk{\Big)}\bU_i^T\Big\}
\right)\vecof(\bSigma-\bU_i\bU_i^T)^T\\[1ex]
&&\qquad\qquad -\big(\bU_i\bU_i^T\big)\otimes\big(\HscrAAi^{-1}\big)\Big|\bX_i\Big].
\end{eqnarray*}
}
Next, note that 
\begin{equation}
\bM_{11}^{-1}=\left[
\begin{array}{cc}
\myOP(m^{-1})\bone_{\dR}^{\otimes2} & \myOP\{(mn)^{-1}\}\bone_{\dR}\bone_{\tcbk{\dB}}^T\\[1ex]
\myOP\{(mn)^{-1}\}\bone_{\dB}\bone_{\tcbk{\dR}}^T & \myOP\{(mn)^{-1}\}\bone_{\dB}^{\otimes2} 
\end{array}
\right].
\label{eq:BarnBluff}
\end{equation}
Given the orders of magnitude in (\ref{eq:WidnesFactory}) and (\ref{eq:BarnBluff}),
from expansion of $\bM_{12}^T\bM_{11}^{-1}\bM_{12}$ 
it is apparent that its dominant $\myOP(m/n)$ contribution is from 
{\setlength\arraycolsep{1pt}
\begin{eqnarray*}
&&\left\{{\displaystyle\sumim} E(\bSBi\bSCi^T|\bX_i)\right\}^T
\left\{{\displaystyle\sumim}E(\bSBi\bSBi^T|\bX_i)\right\}^{-1}
{\displaystyle\sumim} E(\bSBi\bSCi^T|\bX_i)\\[1ex]
&&\qquad\qquad\qquad={\displaystyle\frac{\phi m}{n}}\bD_{\dR}^T(\bSigma^{-1}\otimes\bSigma^{-1})
\bKscBC\bKscBB^{-1}\bKscBC(\bSigma^{-1}\otimes\bSigma^{-1})\bD_{\dR}
+\myoP(m/n)\bone_{\dR(\dR+1)/2}^{\otimes2}
\end{eqnarray*}
}
where 
{\setlength\arraycolsep{1pt}
\begin{eqnarray*}
\bKscBC
&\equiv&\frac{1}{2m}\sumim E\Big[\Big\{\HscrABi^T\HscrAAi^{-1}\bSigma^{-1}\bU_i
-\HscrdAABi\bigstar\HscrAAi^{-1}\\[1ex]
&&\qquad\qquad\qquad+\HscrABi^T\HscrAAi^{-1}
\Big(\HscrdAAAi\bigstar\HscrAAi^{-1}\Big)\Big\}\Big|\bX_i\Big]
\vecof(\bU_i\bU_i^T-\bSigma)^T
\end{eqnarray*}
}
is a matrix with all entries being $\myOP(1)$. Hence, if we let 
$\bAsc\equiv\smhalf\bD_{\dR}^T\big((\bSigmaZero)^{-1}\otimes(\bSigmaZero)^{-1}\big)\bD_{\dR}$ and
$$\bBsc\equiv\bD_{\dR}^T
\big((\bSigmaZero)^{-1}\otimes(\bSigmaZero)^{-1}\big)\big\{\bKscCC^0
+(\bKscBC^0)^T(\bKscBB^0)^{-1}\bKscBC^0\big\}
\big((\bSigmaZero)^{-1}\otimes(\bSigmaZero)^{-1}\big)\bD_{\dR}$$
then 
{\setlength\arraycolsep{1pt}
\begin{eqnarray*}
\Cov\big(\vech(\bSigmaMLE)|\bXsc\big)&=&\frac{1}{m}\left(\bAsc-\frac{\phi}{n}\bBsc\right)^{\tcbk{-1}}
+\myoP\{(mn)^{-1}\}\bone_{\dR(\dR+1)/2}^{\otimes2}\\[1ex]
&=&\frac{1}{m}\bAsc^{-1}+\frac{\phi}{mn}\bAsc^{-1}\bBsc\bAsc^{-1}
+\myoP\{(mn)^{-1}\}\bone_{\dR(\dR+1)/2}^{\otimes2}.
\end{eqnarray*}
}
From (\ref{eq:MandNfive}),
$$\bAsc^{-1}=2\bD_{\dR}^+(\tcbk{\bSigmaZero}\otimes\tcbk{\bSigmaZero})\bD_{\dR}^{+T}.$$
To simplify $\bAsc^{-1}\bBsc\bAsc^{-1}$, we use (\ref{eq:MandNfour}) and \tcbk{\eqref{eq:ACkronBD}}
to obtain
\begin{equation}
{\setlength\arraycolsep{1pt}
\begin{array}{rcl}
\Cov\big(\vech(\bSigmaMLE)|\bXsc\big)
&=&\displaystyle{\frac{2\bD_{\dR}^+(\bSigmaZero\otimes\bSigmaZero)\bD_{\dR}^{+T}}{m}}
+\displaystyle{\frac{4\phi\bD_{\dR}^+
\{\bKscCC^0+(\bKscBC^0)^T(\bKscBB^0)^{-1}\bKscBC^0\}\bD_{\dR}^{+T}}{mn}}\\[1.5ex]
&&\qquad\qquad+\myOP(m^{-1}n^{-2})\bone_{\dR(\dR+1)/2}^{\otimes2}.
\end{array}
}
\label{eq:CovVechNonPop}
\end{equation}

\section{Population Forms of Covariance Matrix Second Terms}

In the previous section, the second terms of the asymptotic covariance matrices 
of $\bbetaMLE$ and $\vech\big(\bSigmaMLE\big)$ are stochastic. However, under
relatively mild moment conditions such as \tcbk{assumption} (A3) of
Jiang \textit{et al.} (2022), these terms converge in probability to
deterministic population forms. In this section we determine these limiting forms. 

A re-writing of the $\bKscAA$ quantity is 
{\setlength\arraycolsep{1pt}
\begin{eqnarray*}
\bKscAA&\equiv&\frac{1}{m}\sumim E\Big[\bU_i\bU_i^T\bSigma^{-1}\big(\smoon\HscrAAi\big)^{-1}         
+\big(\smoon\HscrAAi\big)^{-1}\bSigma^{-1}\bU_i\bU_i^T-\big(\smoon\HscrAAi\big)^{-1}\\[0.5ex]
&&\qquad\qquad+\big(\smoon\HscrAAi\big)^{-1}
\Big\{\big(\smoon\HscrdAAAi\big)\bigstar\big(\smoon\HscrAAi\tcbk{\big)}^{-1}\big)\Big\}\bU_i^T\\[0.5ex]
&&\qquad\qquad+\bU_i\Big\{\big(\smoon\HscrdAAAi\big)\bigstar\big(\smoon\HscrAAi\big)^{-1}\Big\}^T
\big(\smoon\HscrAAi\big)^{-1}
\Big|\bX_i\Big].
\end{eqnarray*}
}
Since
$$E\big(\smoon\HscrAAi|\bX_i\big)\convprob\bOmegaAA(\bU_i)\quad\mbox{and}\quad
E\big(\smoon\HscrdAAAi|\bX_i\big)\convprob\bOmegadAAA(\bU_i)\
$$
we have, under relatively mild conditions (see e.g. Lemma \tcbk{A}1 of Jiang \textit{et al.}, 2022),
{\setlength\arraycolsep{1pt}
\begin{eqnarray*}
\bKscAA&\convprob&\frac{1}{m}\sumim E\Big[\bU_i\bU_i^T\bSigma^{-1}\bOmegaAA(\bU_i)^{-1}         
+\bOmegaAA(\bU_i)^{-1}\bSigma^{-1}\bU_i\bU_i^T-\bOmegaAA(\bU_i)^{-1}\\[0.5ex]
&&\qquad\qquad+\bOmegaAA(\bU_i)^{-1}
\Big\{\bOmegadAAA(\bU_i)\bigstar\bOmegaAA(\bU_i)^{-1}\Big\}\bU_i^T\\
&&\qquad\qquad+\bU_i\Big\{\bOmegadAAA(\bU_i)\bigstar\bOmegaAA(\bU_i)^{-1}\Big\}^T
\bOmegaAA(\bU_i)^{-1}\Big]\\[1ex]
&=&E\Big[\bU\bU^T\bSigma^{-1}\bOmegaAA(\bU)^{-1}         
+\bOmegaAA(\bU)^{-1}\bSigma^{-1}\bU\bU^T-\bOmegaAA(\bU)^{-1}\\[0.5ex]
&&\ \ +\bOmegaAA(\bU)^{-1}
\Big\{\bOmegadAAA(\bU)\bigstar\bOmegaAA(\bU)^{-1}\Big\}\bU^T
+\bU\Big\{\bOmegadAAA(\bU)\bigstar\bOmegaAA(\bU)^{-1}\Big\}^T
\bOmegaAA(\bU)^{-1}\Big]\\[1ex]
&=&\bLambdaAA
\end{eqnarray*}
}
where $\bLambdaAA$ is as defined in Section \ref{sec:theNotation}.
Analogous arguments lead to 
$$\bKscAB\convprob\bLambdaAB,\quad \bKscBB\convprob E\{\MATRIXsix\},\quad \bKscBC\convprob\bPhi
\ \ \mbox{and}\ \ \bD_{\dR}^+\bKscCC\bD_{\dR}^{+T}
\convprob\frac{1}{2}E\big\{\MATRIXnine-2\MATRIXeight\big\}
$$
where $\MATRIXeight$, $\MATRIXnine$ and $\bLambdaAB$ are as defined in Section \ref{sec:theNotation}. 
It follows that the deterministic forms of the order $(mn)^{-1}$ terms match those stated 
in (\ref{eq:mainResult}).

\section{The Gaussian Response Special Case}\label{sec:GaussCase}

For the Gaussian response special case of (\ref{eq:theModel})
the two-term covariance matrix expressions simplify considerably.
The main reason is that, for the Gaussian case, $b''(x)=1$ 
and $b'''(x)=0$. These facts imply that 
$$\bOmegaAA(\bU)=E(\bXA\bXA^T),\quad\bOmegaAB(\bU)=E(\bXA\bXB^T),
\quad \bOmegaBB(\bU)=E(\bXB\bXB^T)$$
and all entries of the three-dimensional arrays $\bOmegadAAA(\bU)$
and $\bOmegadAAB(\bU)$ are exactly zero.

\subsection{The $\Cov\big(\bbetaMLE|\bXsc\big)$ Approximation}

For the Gaussian response situation
$$\bLambdaAA=E(\bXA\bXA^T)^{-1},\quad \bLambdaAB=E(\bXA\bXA^T)^{-1}E(\bXA\bXB^T)
$$
and
$$E\{\MATRIXsix\}=E(\bXB\bXB^T)-E(\bXA\bXB^T)^TE(\bXA\bXA^T)^{-1}E(\bXA\bXB^T).$$
Therefore,
$$
\left[
\begin{array}{cc}
\bLambdaAA^{-1}             &\bLambdaAA^{-1}\bLambdaAB \\[2ex]
\bLambdaAB^T\bLambdaAA^{-1} &\bLambdaAB^T\bLambdaAA^{-1}\bLambdaAB
                              +E\big\{\MATRIXsix\big\}
\end{array}
\right]
=
E\left[
\begin{array}{cc}
\bXA\bXA^T & \bXA\bXB^T\\[1ex] 
\bXB\bXA^T & \bXB\bXB^T 
\end{array}
\right]
=E(\bX\bX^T).
$$
Hence, for the Gaussian special case 
$$\Cov\big(\bbetaMLE|\bXsc\big)=
\frac{1}{m}\left[
\begin{array}{cc}
\bSigmaZero &\quad \bO \\[1ex]
\bO     &\quad  \bO
\end{array}
\right]
+\frac{\phi\big\{E(\bX\bX^T)\big\}^{-1}\{1+\myoP(1)\}}{mn}.
$$
\tcbk{This result generalises the two-term expansion of
  $\Var\big(\betaAMLE|\mathcal{X}\big)$
  provided in Section 3.5 of McCulloch \textit{et al.} (2008)
  for the $\dR=\dB=1$ and $\bXA=1$ special case.}

\subsection{The $\Cov\big(\vech(\bSigmaMLE)|\bXsc\big)$ Approximation}

As shown in, for example, Section 4.3 of Wand (2002) there is exact orthogonality
between $\bbeta$ and $\bSigma$ in the Gaussian case. This means that 
$\bPhi=\bO$ and, hence, the second term of $\Cov\big(\vech(\bSigmaMLE)|\bXsc\big)$ is
\begin{equation}
\frac{2\phi}{mn}
E\Big\{\vech(\bSigma-\bU\bU^T)\VECTORfour^T+\VECTORfour\vech(\bSigma-\bU\bU^T)^T-2\MATRIXeight\Big\}
\label{eq:Bitzberg}
\end{equation}
where $\VECTORfour$ and $\MATRIXeight$ simplify to 
$$\VECTORfour=\bD_{\dR}^+\vecof\Big(\{E(\bXA\bXA^T)\}^{-1}\bSigma^{-1}\big(\bSigma- \bU\bU^T\big)\Big)$$
and
$$\MATRIXeight=\bD_{\dR}^+\big[(\bU\bU^T)\otimes\{E(\bXA\bXA^T)\}^{-1}\big]\bD_{\dR}^{+T}.$$
We immediately have 
$$E\{\MATRIXeight\}=\bD_{\dR}^+\big[\bSigma\otimes\{E(\bXA\bXA^T)\}^{-1}\big]\bD_{\dR}^{+T}.$$
The reduction of the other expectations in (\ref{eq:Bitzberg}) is less immediate and benefits
from Theorem 4.3(iv) of Magnus \myand Neudecker (1979) as well as (\ref{eq:MandNtwo}).
However, such a pathway leads to 
$$E\Big\{\vech(\bSigma-\bU\bU^T)\VECTORfour^T+\VECTORfour\vech(\bSigma-\bU\bU^T)^T\Big\}
=4\bD_{\dR}^+\big[\bSigma\otimes\{E(\bXA\bXA^T)\}^{-1}\big]\bD_{\dR}^{+T}.
$$
On combining the components of (\ref{eq:Bitzberg}) we arrive at 
{\setlength\arraycolsep{1pt}
\begin{eqnarray*}
\Cov\big(\vech(\bSigmaMLE)|\bXsc\big)
&=&\frac{2\bD_{\dR}^+(\bSigmaZero\otimes\bSigmaZero)\bD_{\dR}^{+T}}{m}
+\frac{4\phi\bD_{\dR}^+\big[\bSigmaZero\otimes\{E(\bXA\bXA^T)\}^{-1}\big]\bD_{\dR}^{+T}
\{1+\myoP(1)\}}{mn}.
\end{eqnarray*}
}

\tcbk{\section{Additional Simulation Exercise Figure}\label{sec:additFig}}

\begin{figure}[!ht]
\includegraphics[width=\textwidth]{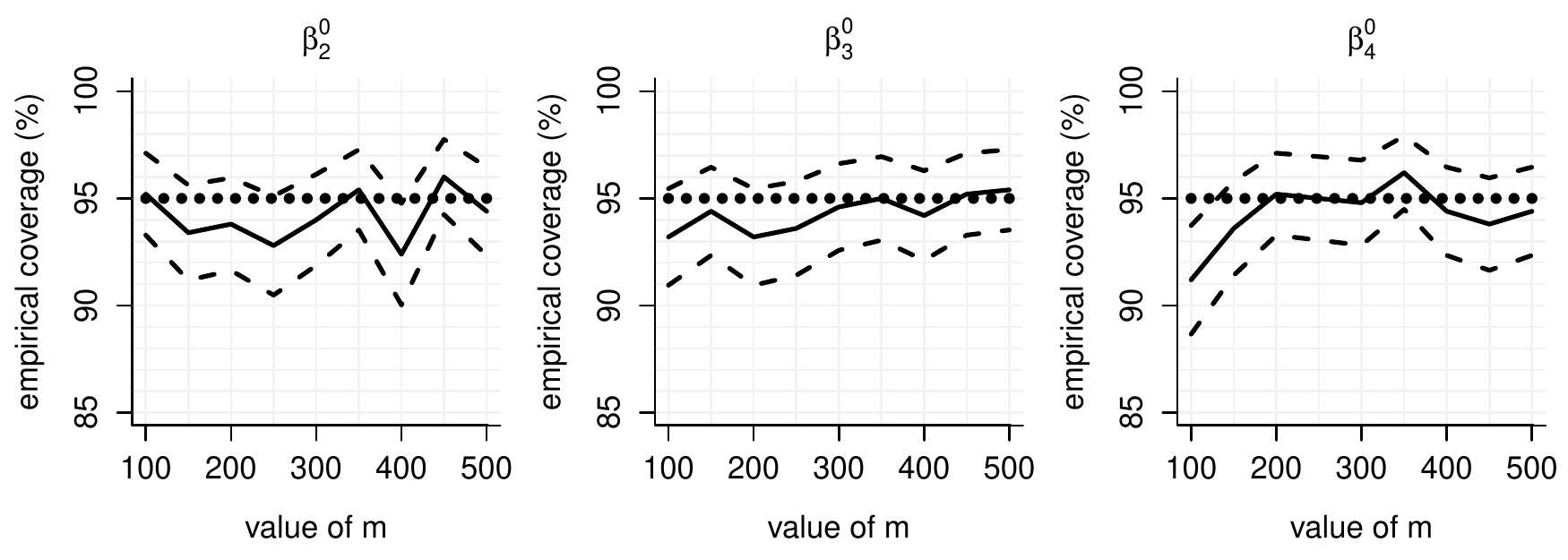}
\caption{\tcbk{\it Empirical coverage of confidence intervals from the simulation
exercise described in Section \ref{sec:numericRes}. Each panel corresponds to a fixed effect model
parameter that is not impacted by second term asymptotic improvements.
The advertised coverage level is fixed at 95\% and is indicated by a
horizontal dotted line in each panel. The solid curves show, dependent on 
the number of groups $m$, the empirical coverage levels for confidence 
intervals that use a one-term asymptotic variance approximation. 
The dashed curves correspond to plus and minus two standard errors of
the sample proportions. The within-group sample size, $n$, is fixed at $m/10$.
}}
\label{fig:MBWcovSuppFig} 
\end{figure}

\tcbk{Figure \ref{fig:MBWcovSuppFig} refers to the simulation exercise described in 
Section \ref{sec:numericRes} and compares the empirical coverages of confidence intervals
with advertised levels of 95\% for the parameters of (\ref{eq:mainResult}) that
are not affected by second term improvement. It is clear from Figure \ref{fig:MBWcovSuppFig} 
that the simple one-term asymptotic variances lead to good coverages for $\betaTwoZero$, 
$\betaThreeZero$ and $\betaFourZero$, even for lower sample size situations.} 

\section*{References}

\bib
Jiang, J., Wand, M.P. \myand Bhaskaran, A. (2022).
Usable and precise asymptotics for generalized linear mixed model analysis
and design. \textit{Journal of the Royal Statistical Society, Series B},
\textbf{84}, 55--82.

\bib
Magnus, J.R. and Neudecker, H. (1979).
The commutation matrix: some properties and applications.
\textit{The Annals of Statistics}, {\bf 7}, 
381--394.

\bib
Magnus, J.R. and Neudecker, H. (1999).
\textit{Matrix Differential Calculus. Revised Edition.}
Chichester, U.K.: John Wiley \& Sons.

\bib
Miyata, Y. (2004).
Fully exponential Laplace approximation using asymptotic modes.
\textit{Journal of the American Statistical Association}, {\bf 99}, 
1037--1049.

\bib
Pace, L. and Salvan, A. (1997). \textit{Principles of Statistical Inference
from a Neo-Fisherian Perspective.} Singapore: World Scientific Publishing
Company.

\bib
Wand, M.P. (2002). Vector differential calculus in statistics.
{\it The American Statistician}, {\bf 56}, 55--62.

\end{document}